\documentclass[12pt]{amsart}
\newif\ifkqboxed      \kqboxedtrue

\makeatletter
\@for\kq@pkg:=lmodern,microtype,mathtools,tikz,xcolor,booktabs,enumitem,hyperref,thmtools,mdframed,cleveref\do{%
  \IfFileExists{\kq@pkg.sty}{}{%
    \PackageError{preamble-ams}{Required package '\kq@pkg' is not installed
      on this machine}{Install it (tlmgr install \kq@pkg) or use a full
      TeX Live / MacTeX; the document does not build without it.}}}
\makeatother

\usepackage[a4paper,margin=3.0cm]{geometry}

\usepackage{lmodern}
\usepackage[T1]{fontenc}

\usepackage{microtype}

\usepackage{amsmath,amssymb,amsthm}
\usepackage{mathtools}

\allowdisplaybreaks[1]

\DeclareMathOperator{\Tr}{Tr}
\DeclareMathOperator{\rk}{rk}      
\newcommand\RG{\mathrm{RG}}

\newcommand\fq{\mathfrak{q}}
\newcommand\fg{\mathfrak{g}}
\newcommand\bZ{\mathbb{Z}}
\newcommand\bR{\mathbb{R}}
\newcommand\bC{\mathbb{C}}

\newcommand\cT{{\mathcal T}}

\newcommand\cM{{\mathcal M}}
\newcommand\IR{\mathrm{IR}}
\newcommand\UV{\mathrm{UV}}
\newcommand\op{\mathrm{op}}
\newcommand\fsl{\mathfrak{sl}}

\newcommand{\cocy}{\mathsf{C}}  
\newcommand{\cocha}{\mathsf{D}}

\newcommand{\cG}{\mathcal{G}} 
\newcommand{\cE}{\mathcal{E}} 

\usepackage{tikz}
\usetikzlibrary{decorations.markings}
\usetikzlibrary{arrows.meta}
\tikzset{
    midarrow/.style={postaction={decorate},decoration={
            markings,
            mark=at position .5 with {\arrow{stealth}}
        }
    }
}

\usepackage{booktabs}
\usepackage{enumitem}

\usepackage{xcolor}
\definecolor{kqrule}{HTML}{1F4E79}    
\definecolor{kqtint}{HTML}{F2F6FA}    
\definecolor{kqdefrule}{HTML}{2E6E4E} 
\definecolor{kqdeftint}{HTML}{F2F8F4} 
\definecolor{kqremrule}{HTML}{8A929E} 
\definecolor{kqlink}{HTML}{1F4E79}

\usepackage{hyperref}                
\hypersetup{
  colorlinks = true,
  linkcolor  = kqlink,
  citecolor  = kqlink,
  urlcolor   = kqlink,
  linktocpage= true,
  bookmarksnumbered = true,
}

\newcommand{\affiliation}[2][]{\address{#2}}
\newcommand{\emailAdd}[1]{\email{#1}}
\makeatletter
\gdef\kq@amsabs{}
\renewcommand{\abstract}[1]{\gdef\kq@amsabs{#1}}
\AtBeginDocument{%
  \ifx\kq@amsabs\@empty\else
    \global\setbox\abstractbox=\vtop{%
      \normalfont\Small
      \list{}{\labelwidth\z@ \leftmargin3pc \rightmargin\leftmargin
        \listparindent\normalparindent \itemindent\z@
        \parsep\z@ \@plus\p@}%
      \item[\hskip\labelsep\scshape\abstractname.]\kq@amsabs
      \endlist}%
  \fi}
\makeatother
\DeclareOldFontCommand{\cal}{\normalfont\rmfamily}{\mathcal}

\theoremstyle{plain}
\newtheorem{theorem}{Theorem}[section]
\makeatletter
\@ifundefined{newcounteralias}{}{%
  \AddToHook{package/thmtools/after}{%
    \@ifpackagelater{thmtools}{2023/05/05}{}{%
      \renewcommand\thmt@autorefsetup{%
        \expandafter\def\csname\thmt@envname autorefname\expandafter\endcsname
          \expandafter{\thmt@thmname}%
      }%
    }%
  }%
}
\makeatother

\usepackage{thmtools}

\declaretheorem[name=Example,     style=definition, sibling=theorem]{example}
\ifkqboxed
  \usepackage[framemethod=TikZ]{mdframed}
  \mdfdefinestyle{kqclaim}{%
    linewidth=0pt,
    leftline=true, rightline=false, topline=false, bottomline=false,
    innerleftmargin=10pt, innerrightmargin=10pt,
    innertopmargin=7pt, innerbottommargin=7pt,
    middlelinewidth=2.5pt, roundcorner=0pt,
    skipabove=\topsep, skipbelow=\topsep,
  }
  \declaretheorem[
    name=Conjecture, style=definition,
    mdframed={style=kqclaim,
              backgroundcolor=kqtint, linecolor=kqrule},
  ]{conjecture}
  \declaretheorem[
    name=Definition, style=definition, sibling=theorem,
    mdframed={style=kqclaim,
              backgroundcolor=kqdeftint, linecolor=kqdefrule},
  ]{definition}
  \declaretheorem[
    name=Remark, style=remark, sibling=theorem,
    mdframed={style=kqclaim,
              linecolor=kqremrule, middlelinewidth=1.5pt},
  ]{remark}
\else
  \declaretheorem[name=Conjecture, style=definition]{conjecture}
  \declaretheorem[name=Definition, style=definition,
                  sibling=theorem]{definition}
  \declaretheorem[name=Remark, style=remark,
                  sibling=theorem]{remark}
\fi

\makeatletter
\providecommand\@dotsep{4.5}
\makeatother

\newcommand{\listofconjectures}{%
  \begingroup
  \renewcommand{\listtheoremname}{List of conjectures}%
  \listoftheorems[ignoreall,show={conjecture},onlynamed={conjecture}]%
  \endgroup}

\makeatletter
\AtBeginDocument{%
  \let\kq@remark@env\remark
  \def\kq@remark@name{remark}%
  \renewcommand\remark{%
    \edef\kq@cur{\@currenvir}%
    \ifx\kq@cur\kq@remark@name
      \expandafter\kq@remark@env
    \else
      \expandafter\kq@remark@bare
    \fi}%
}
\newcommand\kq@remark@bare{%
  \@latex@warning{bare \string\remark\space used as a command; wrap it
    in \string\begin{remark}...\string\end{remark} to get the boxed
    form}%
  \par\addvspace\topsep
  \refstepcounter{theorem}%
  {\itshape Remark \thetheorem.}\enspace\ignorespaces}
\makeatother

\usepackage[capitalise,nameinlink]{cleveref}
  \crefname{conjecture}{Conjecture}{Conjectures}
  \Crefname{conjecture}{Conjecture}{Conjectures}
  \crefname{axiom}{Axiom}{Axioms}
  \Crefname{axiom}{Axiom}{Axioms}
  \crefname{proposition}{Proposition}{Propositions}
  \Crefname{proposition}{Proposition}{Propositions}
  \crefname{lemma}{Lemma}{Lemmas}
  \Crefname{lemma}{Lemma}{Lemmas}
  \crefname{corollary}{Corollary}{Corollaries}
  \Crefname{corollary}{Corollary}{Corollaries}
  \crefname{example}{Example}{Examples}
  \Crefname{example}{Example}{Examples}
  \crefname{definition}{Definition}{Definitions}
  \Crefname{definition}{Definition}{Definitions}
  \crefname{remark}{Remark}{Remarks}
  \Crefname{remark}{Remark}{Remarks}

\title{$K_\fq$-algebras}
\abstract{We introduce the notion of $K_\fq$-algebra, a conjectural axiomatization of the fusion algebra of rotation-equivariant BPS line defects in a 4d ${\cal N}=2$ Supersymmetric Quantum Field Theory. We also introduce the notion of RG flow of $K_\fq$-algebras, a conjectural axiomatization of the action of Seiberg-Witten type RG flows on line defects. These definitions enrich the theory of cluster algebras, $K$-theoretic Coulomb branch algebras, skein algebras and quantum groups.}

\author{Davide Gaiotto}
\affiliation{Perimeter Institute for Theoretical Physics\\ 31 Caroline Street North, Waterloo, ON N2L 2Y5, Canada}
\emailAdd{dgaiotto@perimeterinstitute.ca}
\date{\today}

\begin{document}
\maketitle
\tableofcontents
\listofconjectures      
\newpage
\section{Introduction and Conclusions}
Four-dimensional Supersymmetric Quantum Field Theories with eight supercharges (4d ${\cal N}=2$ SQFTs) give rise to rich mathematical structures \cite{Seiberg:1994rs,Seiberg:1994aj,Donagi:1995cf,Witten:1988ze,Witten:1994cg,Witten:1997sc,Nekrasov:2002qd,Nekrasov:2003rj,Nekrasov:2009rc,Nekrasov:2010ka,Alday:2009aq,Maulik:2012wi,Negut:2015jza,Gaiotto:2009we,Gaiotto:2010okc,Gaiotto:2009hg,Gaiotto:2010be,Gaiotto:2012rg,Kontsevich:2008fj,Kontsevich:2010px,Kontsevich:2013rda,Gaiotto:2024fso,Alim:2011ae,Cecotti:2011rv,Beem:2013sza,ArabiArdehali:2025fad,Gadde:2011uv,Kapustin:2005py,Kapustin:2007wm,Kapustin:2006hi,Oh:2019bgz,Niu:2021jet,Kapustin:2006pk,Nakajima:2015txa,Braverman:2016wma}. Here we are interested in the ``fusion algebra'' of supersymmetric line defects \cite{Kapustin:2007wm,Kapustin:2006hi,Gaiotto:2010be}:
the algebra spanned over $\bZ[\fq,\fq^{-1}]$ by the (K-theory classes of the) line defects, with a formal parameter $\fq$ keeping track of rotation equivariance. 

Depending on the theory in consideration, the fusion algebra matches algebras of great interest: skein algebras and quantum character varieties \cite{Przytycki2006,Turaev1991,Bullock1996,Przytycki1997,Ben-Zvi:2015jua,BenZvi2016,Jordan:2021hop,Kuperberg1997,Sikora2004,Frohman2020,Bonahon2010,Le:2020eds,Le2023,Alekseev:1994pa,Muller2012,Thurston2013}, including quantum groups \cite{schrader2019clusterrealizationuqmathfrakslnquantum,schrader2026ktheoreticcoulombbranchesquiver}, K-theoretic Coulomb branch algebras \cite{Nakajima:2015txa,Braverman:2016wma,Braverman:2016pwk,Nakajima:2017bdt,Braverman:2018lpy,Cautis_2017,Cautis:2018gzt,Finkelberg:2017qpy,Schrader:2019hzd,Cautis:2023jkq,Bullimore:2015lsa} and many cluster algebras \cite{Fomin:2001mwn,Fomin:2002nsg,Berenstein2003,Fomin:2006cre,Fomin2002,Berenstein2004,Fock:2003jdz,Fock:2003xxy,Fock:2008ixh,Fomin:2007rcq,Gross:2014fwa,Lee2013,Davison2016,Qin2015,Berenstein2012}. 

The initial impetus for this work came from the ``Schur quantization'' program \cite{Gaiotto:2024osr}, which employs certain protected indices $I_{a,b}(\fq)\in \bZ(\!(\fq)\!)$ to define an inner product on the fusion algebra which is conjecturally positive-definite for real $-1<\fq<1$. Inspection of concrete examples of the inner products revealed a curious feature: for a natural choice of ``canonical'' linear basis $\{L_a\}$ the indices take the form
\begin{equation}\label{eq:intro-delta}
  I_{a,b}(\fq) \;=\; \Tr\, L_{\rho(a)} L_b \;=\; \Tr\,  L_b L_{\rho^{-1}(a)}\;=\; \delta_{a,b} + O(\fq)\,,
\end{equation}
where $\Tr$ is a twisted trace and $\rho$ an automorphism of the algebra, both axiomatized below. This is a $\fq \to 0$ version of the positive definiteness statement. 

As we explored the implications of this property on the fusion algebras arising from Coulomb branches, cluster algebras and skein algebras, we encountered several other natural constraints on the $\fq \to 0$ behaviour of various quantities of interest. This paper formalizes such constraints as a collection of abstract mathematical conjectures. As a comparison, the scope of these conjectures seems to be comparable to the scope of the theory of cluster algebras. Via Schur quantization, we also expect them to have interesting implications for the Langlands program. 

We first encountered some of the properties conjectured here as consequences of super-conformal symmetry, when available. Direct inspection of many examples led us to extend the conjectures to non-conformal systems and to propose further constraints based on the overall self-consistency of the problem. As we currently lack physical justifications for the resulting set of conjectures, we segregate the tentative physical arguments to Appendix \ref{app:physics}. 

In conclusion, in this paper we enhance the notion of fusion algebra by formulating a collection of axioms defining a ``$K_\fq$-algebra'' (Definition~\ref{def:kq}). We then describe how the relevant skein algebras, quantum character varieties, K-theoretic Coulomb branch algebras, quantum groups and many cluster algebras can be equipped with the structure of $K_\fq$-algebras. 
We also introduce the notion of an ``RG flow of $K_\fq$-algebras'' (Definition~\ref{def:rg}), which formalizes the effect of Coulomb branch RG flows on fusion algebras. We conjecture that the canonical maps from a physically-motivated quantum cluster variety to any of its quantum torus charts can be equipped with the structure of an RG flow of $K_\fq$-algebras (Conjecture~\ref{conj:upper}). The same is true about certain well-known maps between K-theoretic Coulomb branch algebras, as well as many new examples (Conjecture~\ref{conj:factored}). We further sharpen the 
relation to cluster algebras by a new definition of the Donaldson-Thomas element characterizing the RG flow (Conjecture~\ref{conj:cluster}). 

Upon holomorphic-topological (HT) twist of the SQFT, $K_\fq[\cT]$ can be defined as the rotation-equivariant Grothendieck group of an HT-factorization category $\mathrm{Lines}_{\bC^*}[\cT]$ of topological line defects. 
A well studied example of this structure is the mathematical construction \cite{Cautis:2018gzt,Cautis:2023jkq} of a candidate $\mathrm{Lines}_{\bC^*}[\cT]$ for gauge theories, whose Grothendieck group is the K-theoretic Coulomb branch algebra and conjecturally matches the fusion algebra $K_\fq[G,T^*N]$ of a conventional gauge theory with cotangent-type matter. The properties we identify appear to be naturally associated to a generalization of the HT framework where one explores the (twistor) space of HT structures on $\bR^3$.  More concretely, they suggest a remarkable triangularity structure for the ``bar-invariant simple objects'' in the category $\mathrm{Lines}_{\bC^*}[\cT]$.

The paper is organized as follows: Section~\ref{sec:kq}, the $K_\fq$-algebra axioms, their flavoured variant and the three examples; Section~\ref{sec:coulomb}, K-theoretic Coulomb branch algebras; Section~\ref{sec:rg}, RG flows, Seiberg--Witten flows onto a quantum torus and more general flows derived from it; Section~\ref{sec:skein}, a tentative proposal for skein algebras, where $\rho$ and the trace remain conjectural. The appendices collect the physical motivations (Appendix~\ref{app:physics}) and $K_\fq$-algebras of finite type (Appendix~\ref{app:finite}). All conjectures are listed after the table of contents.

\subsection*{Statement of AI use}
In order to collect supporting evidence for our conjectures, we decided to code a variety of different $K_\fq$-algebras and RG flows between them. This aspect of the project seemed to be a natural testing ground for the powerful agentic AI coding tools which recently became available to the public. An important challenge was to ensure a faithful encoding of the mathematical axioms we proposed,
but also to develop algorithms which could uniquely pin down the previously unknown $K_\fq$-algebra structures possessed by cluster algebras, K-theoretic Coulomb branches and skein algebras. 
The final product of these coding efforts is a public GitHub repository intended as a companion to this paper.\footnote{\url{https://github.com/DGaiotto/KAlgebra}}

The repository contains a reasonably complete documentation of the algorithms we employed and of the large number of $K_\fq$-algebras and RG flows we tested, so we will not attempt to describe all of these in detail in this paper. 
The initial coding work was done with Claude Opus $4.7$ and $4.8$ and required a considerable and detailed involvement on the author's part. The second half of the coding project was then completed with greater ease by Claude Fable 5, Opus 5. Opus 5.5 was used in the final testing stage. 

Within the body of the paper and Appendix A, AI use was restricted to the output of some equations and figures and minor proof-reading, mostly to ensure that the conventions and results matched the 
code output. Appendix B incorporates (with author's rewriting) an AI-generated description of the $K_\fq$-algebra combinatorics, matching the code employed in the companion repository. AI also helped devising the latex preamble for this file. 

The companion repository is entirely AI-written, including the documentation. The author designed the high level object-oriented structure implementing the mathematical framework 
and the associated computational algorithms. The AI handled the actual Python implementation and extensive testing.  
\newpage
\section{$K_\fq$-algebras.} \label{sec:kq}
\begin{definition}[$K_\fq$-algebra]
\label{def:kq}
A \emph{$K_\fq$-algebra} is an algebra $A_\fq$ over $\bZ[\fq, \fq^{-1}]$, free
as a $\bZ[\fq, \fq^{-1}]$-module, endowed with:
\begin{enumerate}[label=(K\arabic*), ref=K\arabic*, leftmargin=*, itemsep=5pt]
  \item\label{ax:bar}
    A bar involution, i.e.\ an antimultiplicative $\bZ$-algebra involution
    sending $\fq \to \fq^{-1}$, $A_\fq \simeq A^\op_{\fq^{-1}}$.
  \item\label{ax:basis}
    A ``canonical'' linear basis $L_a$ over $\bZ[\fq, \fq^{-1}]$ consisting of
    bar-invariant elements and including the identity $1$.
  \item\label{ax:rho}
    An algebra automorphism $\rho$ acting as a permutation
    $L_a \to L_{\rho(a)}$ of the canonical basis elements.
  \item\label{ax:trace}
    A $\rho^2$-twisted trace $\Tr : A_\fq \to \bZ(\!(\fq)\!)$, linear over
    $\bZ[\fq, \fq^{-1}]$, such that
    \begin{equation}\label{eq:orthonormality}
      I_{a,b} \;\equiv\; \Tr L_{\rho(a)} L_b \;=\; \Tr L_{b} L_{\rho^{-1}(a)}
      \;=\; \delta_{a,b} + O(\fq)\,.
    \end{equation}
   \item\label{ax:rhotr}
   The trace satisfies $\Tr L_{\rho(a)} = \Tr L_a$, implying $I_{b,a} = I_{a,b}$. 
\end{enumerate}
\end{definition}

\begin{remark} The axioms allow one to rescale the trace by an overall factor of the
form $1 + O(\fq)$. In physics examples we have a canonical normalization for the
trace. \end{remark}

\begin{remark} The notion of canonical basis and in particular the last entry of axiom \ref{eq:orthonormality} appear with a similar role in \cite{Webster_2014}, up to $\fq_{\mathrm{here}} = q^{-1}_{\mathrm{there}}$.
It would be interesting to pursue the analogy further, especially for $K_\fq$-algebras which are also quantum groups, see Section \ref{sec:uqsl2}. \end{remark}

\begin{remark} \label{rem:maximal} As indicated in the introduction, the notion of $K_\fq$-algebra axiomatized the
expected properties of fusion algebras $K_\fq[\cT]$ of BPS line defects in 4d
${\cal N}=2$ SQFTs $\cT$. The key ``$\delta_{a,b}$'' axiom~\eqref{eq:orthonormality} is expected to hold
only after one maximally extends the collection of line defects available to the
theory, in particular by including direct summands of available line defects.
\end{remark}

\begin{example}[the quantum torus]
\label{ex:qt}
The canonical example of a $K_\fq$-algebra is the quantum torus $Q_\fq(\Gamma)$
associated to a lattice $\Gamma$ equipped with a non-degenerate integral
antisymmetric pairing $\langle \cdot, \cdot \rangle : \Gamma \times \Gamma \to
\bZ$, with generators $X_\gamma$ for $\gamma \in \Gamma$ and relations
\begin{equation}
  X_\gamma X_{\gamma'} = \fq^{\langle \gamma, \gamma' \rangle} X_{\gamma + \gamma'}\,.
\end{equation}
The bar involution fixes the canonical basis elements $X_\gamma$. The
automorphism $\rho$ acts as $\rho(\gamma) = - \gamma$ and squares to the
identity. The trace is
\begin{equation}
  \Tr X_\gamma = (\fq^2)_\infty^{\rk \Gamma}\, \delta_{\gamma,0}
\end{equation}
with a physics motivated overall factor. In particular,
\begin{equation}
  I_{\gamma, \gamma'} = (\fq^2)_\infty^{\rk \Gamma}\, \delta_{\gamma, \gamma'}\,.
\end{equation}
This $K_\fq$-algebra encodes the fusion algebra of BPS line defects in a 4d ${\cal N}=2$ Abelian gauge theory \cite{Gaiotto:2010be}. 
\end{example}

\subsection{The Pentagon algebra, aka $K_\fq([A_1,A_2])$.}
\label{sec:pentagon}
A canonical example of $K_\fq$-algebra $K_\fq([A_1,A_2])$ arises from the $[A_1,A_2]$ Argyres-Douglas theory \cite{Argyres:1995jj,Argyres:1995xn,Gaiotto:2010be,Cecotti:2010fi,Alim:2011ae}. 
\begin{definition}[Pentagon Algebra] \,
\label{def:penta}
\begin{itemize}
	\item $K_\fq([A_1,A_2])$ is generated multiplicatively by five canonical basis elements $L_i$,  
	with relations 
	\begin{align}
		&L_{i+1} L_i \,\,= \fq^2 L_i L_{i+1} \cr
		&L_{i+1} L_{i-1} = 1+\fq L_i \cr
		&L_{i-1} L_{i+1} = 1+\fq^{-1} L_i \, .
	\end{align}
	We use a cyclic notation $L_{i+5}= L_{i}$.  
	\item Canonical basis elements take the form $L_{i;a,b} = \fq^{ab} L_i^a L_{i+1}^b$. 
	\item The automorphism $\rho(L_i) = L_{i+2}$ generates the $\bZ_5$ automorphism group of $K_\fq([A_1,A_2])$.
	\item The trace is determined by $\rho^2$-cyclicity from the seeds:
	 \begin{align}
		\Tr\, 1 &\;=\;  \sum_{n=0}^{\infty} \frac{\fq^{2n(n+1)}}{(\fq^2;\fq^2)_n} \;=\; 1 + \fq^4 + \fq^6 + \fq^8 + \fq^{10} + \ldots \cr
		\Tr\, L_i &\;=\;  -\sum_{n=0}^{\infty} \frac{\fq^{2(n+1)^2-1}}{(\fq^2;\fq^2)_{n}} \;=\;-\fq - \fq^7 -\fq^9 - \fq^{11}  + \ldots \, ,
	\end{align}
	as explained below in greater detail.
\end{itemize}
\end{definition} 
\begin{remark} It is not obvious from this definition that $K_\fq([A_1,A_2])$ exists and satisfies all the $K_\fq$-algebra axioms. The existence as an algebra follows from a well-known presentation as a cluster algebra, given a physical interpretation in \cite{Gaiotto:2010be} in a manner which highlights the roles of the canonical basis and $\rho$. This can be enhanced by the IR formulae for the Schur index \cite{Cordova:2016uwk}, formalized in Definition \ref{def:sw}, to demonstrate the existence of a trace satisfying the remaining axioms. \end{remark}

As an exercise, we will now demonstrate how the trace is fully constrained by the $K_\fq$-algebra axioms, including the key ratio $\Tr L_1/\Tr 1$.  Note the trace is $\rho^2$ and thus $\bZ_5$ invariant. The $\rho^2$-twisted trace relations allow some immediate simplifications: 
\begin{equation} \label{eq:izero}
	 \Tr \,L_i^a L_{i+1}^b =   \Tr L_i^{a+1} L_{i+1}^{b-1} 
\end{equation}
immediately tells us $\Tr\, L_{i;a,b} = \fq^{ab} \Tr L_i^{a+b}$. Then 
\begin{align}
	\Tr\, L^{a+1}_i &= \Tr\, L^a_i L_{i+1} = \fq^{-2a} \Tr\, L_{i+1} L_i^{a} = \fq^{-2a} \Tr \,L_i^{a} L_{i+2} =\cr &= \fq^{-2a} \Tr \,L_i^{a-1}+\fq^{-2a-1} \Tr \,L_i^{a-1} L_{i+1}= \cr &=\fq^{-2a} \Tr \,L_i^{a-1}+\fq^{-2a-1} \Tr \,L_i^{a}\, ,
\end{align}
a beautiful recursion relation which determines all $\rho^2$-twisted traces as linear combinations $\Tr 1$ and $\Tr L_1$ over $\bZ[\fq, \fq^{-1}]$. 

We can argue that $\rho^2$-cyclicity forces no further relations between $\Tr 1$ and $\Tr L_1$. Twisted cyclicity only needs to be tested on generators, and by $\bZ_5$ symmetry it reduces to 
\begin{equation}
	\Tr L_0 L_{i;a,b} = \Tr L_{i;a,b} L_1 \, .
\end{equation}
We thus have five families of checks. The $i=0$ case is equation (\ref{eq:izero}). We find 
\begin{align}
	\Tr L_0 L_{1;a,b} &= \fq^{ab-2a} \Tr L_{1}^a L_0 L_2^b = \fq^{ab-2a} \Tr L_{1}^a L_2^{b-1}+ \fq^{ab-2a-1} \Tr L_{1}^{a+1} L_2^{b-1}= \cr &=\fq^{ab+2b} \Tr L_{1}^{a+1} L_2^{b} =  \Tr L_{1;a,b} L_1 \cr
	\Tr L_0 L_{2;a,b} &= \fq^{ab} \Tr L_0 L_{2}^a L_3^b =\fq^{ab} \Tr L_{2}^{a-1} L_3^b +\fq^{ab-1} \Tr L_1 L_{2}^{a-1} L_3^b 
	= \cr &=\fq^{ab} \Tr L_{2}^{a-1} L_3^b +\fq^{ab-2a+1} \Tr L_{2}^{a-1}  L_1 L_3^b= \cr 
	&=\fq^{ab} \Tr L_{2}^{a-1} L_3^b +\fq^{ab-2a+1} \Tr L_{2}^{a-1}  L_3^{b-1}+\fq^{ab-2a} \Tr L_{2}^{a} L_3^{b-1}  =\cr
	&=\fq^{ab} \Tr L_{2}^{a} L_3^{b-1} +\fq^{ab+2b-1} \Tr L_{2}^{a+1} L_3^{b-1} =\cr
	&=\fq^{ab} \Tr L_{2}^{a} L_3^{b}L_1 =  \Tr L_{2;a,b} L_1 \cr
	\Tr L_0 L_{3;a,b} &= \fq^{ab} \Tr L_0 L_{3}^a L_4^b = \fq^{ab} \Tr L_{3}^{a-1} L_4^b + \fq^{ab+2a-1} \Tr L_{3}^{a-1} L_4^{b+1} = \cr &= \fq^{ab} \Tr L_{3}^{a-1} L_4^b + \fq^{ab-2b+1} \Tr L_{3}^{a-1} L_4^{b-1}+ \fq^{ab-2b} \Tr L_{3}^{a-1} L_4^{b} = \cr
	&= \fq^{ab} \Tr L_{3}^{a-1} L_4^b + \fq^{ab-2b+1} \Tr L_{3}^{a} L_0 L_4^{b-1} =\cr 
	&= \fq^{ab} \Tr L_{3}^{a} L_4^{b-1} + \fq^{ab-1} \Tr L_{3}^{a}  L_4^{b-1} L_0 =  \Tr L_{3;a,b} L_1 \cr
	\Tr L_0 L_{4;a,b} &= \fq^{ab} \Tr L_0 L_{4}^a L_0^b = \fq^{ab+2a } \Tr L_{4}^a L_0^{b+1} = \cr
	 &= \fq^{ab-2b } \Tr L_{4}^{a-1} L_0^{b}+ \fq^{ab-2b-1 } \Tr L_{4}^{a-1} L_0^{b+1}= \fq^{ab-2b } \Tr L_{4}^{a} L_1 L_0^{b}= \cr &=\Tr L_{4;a,b} L_1
\end{align}
The negative powers of $\fq$ in the recursion equation have some interesting implications when combined with the $O(\fq)$ axioms. 
For example,
\begin{equation}
	\Tr \, \rho(L_{i-2}) L_i = \Tr \, L_i^2 = \Tr \,L_{i-1} L_i =  \fq^{-2} \Tr \,L_{i-2} L_i =  \fq^{-2} \Tr\, 1+ \fq^{-3} \Tr \,L_{i} \,
\end{equation}
together with 
\begin{equation}
	\Tr \, \rho(L_{i-2}) L_i  = O(\fq) \,
\end{equation} 
imposes
\begin{equation}
	\Tr \,L_1 = \left(- \fq + O(\fq^4) \right) \Tr \,1\, .
\end{equation} 

More generally, the recursion relation 
\begin{equation}\label{eq:pentarec}
	\Tr \,L_i^n =\fq^{1-2n} \Tr\, L^{n-1}_i + \fq^{2-2n} \Tr\,L^{n-2}_i\,
\end{equation}
leads to expressions
\begin{equation}\label{eq:ABdef}
  \fq^{\,n^2-1}\, \Tr L_i^{\,n} \;=\; A_n(\fq)\, \Tr 1 \;+\; B_n(\fq)\, \Tr L_1 \,,
\end{equation}
which determine the ratio of  $\Tr L_1$ and $\Tr 1$ to order $O(\fq^{n^2})$. In the limit $n \to \infty$ we obtain a remarkable 
continued fraction formula
\begin{equation}
	1 \;-\; \fq\,\frac{\Tr\,L_{i}}{\Tr\,1}
\;=\;
1 + \cfrac{\fq^{2}}{1 + \cfrac{\fq^{4}}{1 + \cfrac{\fq^{6}}{1 + \ddots}}}
\end{equation}

One may recognize the recursion relation for $\Tr L_i^a$ as being very close to the one employed by Schur to study the Yang-Lee $\cM(2,5)$ CFT characters:
\begin{align}
\chi_0(q) &\;=\; \prod_{n=1}^{\infty} \frac{1}{(1-q^{5n-2})(1-q^{5n-3})}
         \;=\; \sum_{n=0}^{\infty} \frac{q^{n(n+1)}}{(q;q)_n}
         \;=\; \cr &\;=\; 1 + q^2 + q^3 + q^4 + q^5 + 2q^6 + 2q^7 + 3q^8 + 3q^9 + 4q^{10} + \ldots \cr
\chi_1(q) &\;=\; \prod_{n=1}^{\infty} \frac{1}{(1-q^{5n-1})(1-q^{5n-4})}
         \;=\; \sum_{n=0}^{\infty} \frac{q^{n^2}}{(q;q)_n}
         \;=\; \cr &\;=\; 1 + q + q^2 + q^3 + 2q^4 + 2q^5 + 3q^6 + 3q^7 + 4q^8 + 5q^9 + \ldots \, .
\end{align}
 and the Rogers-Ramanujan identity they satisfy. Indeed, there are known relations
\begin{align}
\Tr\, 1 &\;=\;  \chi_0(\fq^2)  \cr
\Tr\, L_i &\;=\; \fq^{-1} \chi_0(\fq^2) - \fq^{-1} \chi_1(\fq^2) \, .
\end{align}

The polynomial coefficient $A_n(\fq)$ and $B_n(\fq)$ themselves stabilize as power series as $n \to \infty$, see e.g. Table \ref{tab:pentagon-stabilization}. \begin{table}[t]
  \centering
  \renewcommand{\arraystretch}{1.35}
  \caption{The coefficients of \eqref{eq:ABdef}.  Each row extends the one
    above it; the last row is the limit.}
  \label{tab:pentagon-stabilization}
  \begin{tabular}{@{}r@{\qquad}l@{\qquad}l@{}}
    \toprule
    $n$ & $A_n$ & $B_n$ \\
    \midrule
    $3$ & $\fq$                                   & $1+\fq^{4}$ \\
    $4$ & $\fq+\fq^{7}$                           & $1+\fq^{4}+\fq^{6}$ \\
    $5$ & $\fq+\fq^{7}+\fq^{9}$                   & $1+\fq^{4}+\fq^{6}+\fq^{8}+\fq^{12}$ \\
    $6$ & $\fq+\fq^{7}+\fq^{9}+\fq^{11}+\fq^{17}$ & $1+\fq^{4}+\fq^{6}+\fq^{8}+\fq^{10}+\fq^{12}+\fq^{14}+\fq^{16}$ \\
    \midrule
    $\infty$ & $-\Tr L_1$                         & $\Tr 1$ \\
    \bottomrule
  \end{tabular}
\end{table}
\begin{remark}
\label{rem:miracle}
A curious observation we cannot motivate physically is that $A_\infty = -\Tr L_1$ and $B_\infty = \Tr 1$ respectively, 
so the recursion somehow ``knows'' the traces themselves even though it only constrains their ratio. We will encounter analogous unexpected identities in several other $K_\fq$-algebras examples.
\end{remark}
We refer to Appendix \ref{app:a1a2k} for a detailed description of a larger family $K_\fq([A_1,A_{2k}])$ of (conjectural) $K_\fq$-algebras, which generalize the pentagon algebra and have analogous relations to the $\cM(2,2k+3)$ characters.

\subsection{Flavoured $K_\fq$-algebras.}
\label{sec:flavoured}
We now propose a definition of $K_\fq$-algebras ``flavoured'' by a connected reductive group $G_f$. Denote as $R_{G_f}$ the representation ring of $G_f$ and as $\chi_r$ the finite-dimensional irreducible representations.

\begin{definition}[flavoured $K_\fq$ algebras] \label{def:kq-flavoured}
A \emph{$K_\fq$-algebra flavoured by $G_f$} is an algebra $A_\fq[G_f]$ over $\bZ[\fq, \fq^{-1}] \otimes R_{G_f}$, free
as an $\bZ[\fq, \fq^{-1}]\otimes R_{G_f}$-module but not canonically over the $R_{G_f}$ factor, endowed with 
\begin{enumerate}[label=(F\arabic*), ref=F\arabic*, leftmargin=*, itemsep=5pt]
  \item\label{ax:f-bar}
    A bar involution, i.e.\ an antimultiplicative $\bZ$-algebra involution
    sending $\fq \to \fq^{-1}$, $A_\fq[G_f] \simeq A^\op_{\fq^{-1}}[G_f]$.
  \item\label{ax:f-basis}
    A bar-invariant canonical basis defined over $\bZ[\fq, \fq^{-1}]$, from which one can select, non-canonically, a basis over $\bZ[\fq, \fq^{-1}] \otimes R_{G_f}$. In particular, it includes the basis of characters $\chi_r \in R_{G_f}$ of finite-dimensional $G_f$ irreps $r$. The ambiguity in selecting the basis is associated to tensor products with invertible (1d) irreps of $G_f$. 
  \item\label{ax:f-rho}
    An algebra automorphism $\rho$ which is twisted linear in $R_{G_f}$, with
    $\rho(\chi_r) = \chi_{r^\vee}$, and acts as a permutation of the canonical basis elements.
  \item\label{ax:f-trace}
    A $\rho^2$-twisted trace $\Tr : A_\fq[G_f] \to \bZ(\!(\fq)\!) \otimes R_{G_f}$
    linear over $\bZ[\fq, \fq^{-1}] \otimes R_{G_f}$.  Defining 
    \begin{equation}\label{eq:f-index}
      I_{a,b} \;\equiv\; \Tr L_{\rho(a)} L_b \;=\; \Tr L_{b} L_{\rho^{-1}(a)} \,,
    \end{equation}
    we have
    \begin{equation}\label{eq:f-integrality}
      I_{a,b} \;\in\; R_{G_f} + \fq\, R_{G_f}[[\fq]] \,.
    \end{equation}
    Denoting as $I^{(1)}_{a,b}$ the identity summand, we also have
    \begin{equation}\label{eq:f-orthonormality}
      I^{(1)}_{a,b} = \delta_{a,b} + O(\fq)\,.
    \end{equation}
     \item\label{ax:frhotr}
   The trace satisfies $\Tr L_{\rho(a)} = \rho(\Tr L_a)$, implying $I_{b,a} = \rho(I_{a,b})$. The outer $\rho$ just acts on the flavour ring. 
  \item\label{ax:f-forget}
    We have a forgetful map $A_\fq[G_f] \to A_\fq$ to a $K_\fq$-algebra $A_\fq$, mapping a representation to its dimension. The forgetful map identifies the basis of $A_\fq[G_f]$ over $\bZ[\fq, \fq^{-1}] \otimes R_{G_f}$ with the canonical basis of $A_\fq$, and commutes with $\rho$ and with the trace.
\end{enumerate}
\end{definition}
This definition formalizes the fusion algebra $K_\fq[\cT;G_f]$ of BPS line defects in a 4d ${\cal N}=2$ theory $\cT$ with a flavour symmetry $G_f$. In a categorical language, $K_\fq[\cT;G_f]$ should be the $\bC^* \times G_f$-equivariant Grothendieck group of line defects. The forgetful map 
recovers the conventional Grothendieck group $K_\fq[\cT]$. 

\begin{remark}
\label{rem:flavour-change}
	Given a connected reductive group $H_f$ and any group homomorphism $H_f \to G_f$, the map $R_{G_f} \to R_{H_f}$ allows one to define a $K_\fq$-algebra $A_\fq[H_f]$ flavoured by $H_f$, with a map from $A_\fq[G_f]$ generalizing the forgetful map.
\end{remark}

\begin{example}[flavoured quantum torus]
\label{ex:qt-flavoured}
The quantum torus $Q_\fq(\Gamma)$ associated to a lattice $\Gamma$ with pairing
kernel $\Gamma_f$ is a canonical example of $K_\fq$-algebra flavoured by the
algebraic torus $T_{\Gamma_f}$, identifying $X_{\gamma_f}$ with characters of
1d irreps.  The flavoured trace becomes
\begin{equation}\label{eq:qt-flavoured-trace}
  \Tr X_\gamma = (\fq^2)_\infty^{\rk(\Gamma/\Gamma_f)}\, \chi_{\gamma}\,
                 \delta_{\gamma \in \Gamma_f}\,.
\end{equation}
\end{example}

\begin{remark}
\label{rem:flavour-comments}
A potential obstruction to our setup is the possibility that certain line defects may break the $G_f$ symmetry, i.e. not be $G_f$-equivariant. In the current context of 4d ${\cal N}=2$ half-BPS line defects, we are not aware of any examples 
of such phenomenon, at least for connected $G_f$. We will thus assume that every such line defect is automatically equivariant and that canonical basis elements in $K_\fq[\cT]$ admit a lift to $K_\fq[\cT;G_f]$. The lift may not be completely canonical, but different lifts will be related by tensoring with a one-dimensional  
representation of $G_f$ as indicated above. In order for our formalization to hold, it may be necessary to extend the available collection of line defects in the theory, 
as well as replace $G_f$ with a finite cover. We will assume that this is always possible for the theories we consider. The case of disconnected $G_f$ is certainly interesting: for example, it should be possible to discuss a $\bZ_5$-equivariant version of the Pentagon algebra trace. We leave that to future work. 
\end{remark}

\subsection{The $U_\fq(\fsl_2)$ quantum group as a flavoured $K_\fq$-algebra}
\label{sec:uqsl2}

The next example is an $SU(2)$-flavoured $K_\fq$-algebra which arises from
an $U(1)$ gauge theory with $T^* \bC^2$ matter, i.e. $\mathrm{SQED}_2$. 
It is isomorphic, as an algebra, to the central quotient of
$U_\fq(\fsl_2)$. 
\begin{definition}[$U_\fq(\fsl_2)$ $K_\fq$-algebra] \,
\begin{itemize}
	\item The algebra is multiplicatively generated by the Cartan $K^{\pm 1}$ and the
	canonical basis elements $E,F$, with relations
	\begin{align}
		K E &= \fq^{-2} E K \,,                    & K F &= \fq^{2} F K \,, \\
		E F &= \chi_1 + \fq K + \fq^{-1} K^{-1}\,, & F E &= \chi_1 + \fq^{-1} K + \fq K^{-1} \,,
	\end{align}
	so that $[E,F] = (\fq-\fq^{-1})(K-K^{-1})$ and the quadratic Casimir is the fundamental $SU(2)$
	flavour character $\chi_1$ itself. 
	\item Other canonical basis elements take the form $E_{a,b} = \fq^{-ab} E^a K^b$, $F_{a,b}=\fq^{ab} F^a K^b$ and the Cartan powers $K^n$. These powers of $\fq$ are fixed by bar-invariance.
	\item We have 
	\begin{equation}
		\rho(K)=K^{-1} \,, \qquad \rho(E)=\fq^{-1}FK^{-1} \,, \qquad \rho(F)=\fq\,KE \,,
	\end{equation}
	with $\rho$ fixing the characters $\chi_k$ ($SU(2)$ representations are self-dual). 
	\item The trace vanishes on $E_{a,b}$ and $F_{a,b}$.  Traces of $K^n$ are collected by the
	generating function
	\begin{equation}\label{eq:uqsl2trace}
		G(x,\mu) \equiv \sum_{n\in\bZ} \Tr K^n \; x^n = (\fq^2;\fq^2)_\infty^2 \;
		E_\fq(\mu x)\,E_\fq(\mu^{-1} x)\,E_\fq(\mu x^{-1})\,E_\fq(\mu^{-1} x^{-1}) \, ,
	\end{equation}
	where $\chi_1 = \mu+ \mu^{-1}$ and $E_\fq(x) = (-\fq x;\fq^2)_\infty^{-1}$. This formula is a special case of the Schur index for general gauge theories, which we will review momentarily. 
\end{itemize}
\end{definition}
Unlike the pentagon algebra $\rho$, here $\rho$ is of infinite order: it is Lusztig's braid symmetry,
covering the order-two Weyl reflection $K\mapsto K^{-1}$ but acting as a braid
group $\bZ$. \footnote{Our general terminology differs from $U_\fq(\fsl_2)$ conventions in important ways. The bar involution ~\ref{ax:bar} is antimultiplicative and fixes $K$, whereas the conventional bar involution of $U_\fq(\fsl_2)$ is an algebra automorphism and sends $K \mapsto K^{-1}$. Our canonical basis \ref{ax:basis} is not the canonical basis of Lusztig and Kashiwara.}
	
It is instructive to see again how the trace is constrained by $\rho^2$-twisted cyclicity. Applying it to $\Tr(K^n E F)$
in the form $\Tr(K^n E F)=\Tr(\rho^2(F)\,K^n E)$, and using $EF=\chi_1+\fq K+\fq^{-1}K^{-1}$ together with the $K$-weights of $E,F$,
gives the recursion
\begin{equation}\label{eq:uqsl2rec}
	\fq^{2n}\big(\Tr K^n + \fq\,\chi_1\Tr K^{n+1} + \fq^{2}\,\Tr K^{n+2}\big)
	= \Tr K^n + \fq\,\chi_1\Tr K^{n-1} + \fq^{2}\,\Tr K^{n-2} \, .
\end{equation}
i.e. 
\begin{align}\label{eq:uqsl2con}
	(1-\fq^{2n}) \Tr K^n &= \fq^{2n+1}\chi_1\Tr K^{n+1} + \fq^{2n+2}\Tr K^{n+2} -  \fq\chi_1\Tr K^{n-1} - \fq^{2}\Tr K^{n-2} \qquad n>0 \cr
	(1-\fq^{-2n}) \Tr K^n &= -\fq \chi_1\Tr K^{n+1} - \fq^{2}\Tr K^{n+2} +  \fq^{1-2n} \chi_1\Tr K^{n-1} + \fq^{2-2n}\Tr K^{n-2} \qquad n<0  \, .
\end{align}
Given any $\Tr 1 = 1 + O(\fq)$, we can iterate this recursion from the $O(\fq)$ approximation $\Tr K^n = O(\fq)$ to get an unique solution at order $O(\fq^2)$, then $O(\fq^3)$, etc. It is easy to verify $G(x,\mu)$ produces such a solution. 

We have thus proven that the central quotient of $U_\fq(\fsl_2)$ can be equipped with the structure of an $SU(2)$-flavoured $K_\fq$-algebra.

\begin{remark} Due to the relation between $K_\fq$ algebras and quantum character varieties partially reviewed in Section \ref{sec:skein}, and the relation between 
quantum character varieties and quantum groups described e.g. in \cite{Gaiotto:2014lma,schrader2019clusterrealizationuqmathfrakslnquantum} for type A, we expect some variant of $U_\fq(\fg)$ 
to admit the structure of a $K_\fq$ algebra for every ADE $\fg$ and perhaps every reductive $\fg$.  \end{remark}
\subsection{The $K_\fq([A_1,D_3])$ algebra.}
\label{sec:a1d3}
Our next example is another $SU(2)$-flavoured $K_\fq$-algebra, arising from the
$[A_1,D_3]$ Argyres-Douglas theory. Like the pentagon algebra it enjoys a cyclic symmetry generated by $\rho$, now $\bZ_3$.
\begin{definition}[The punctured triangle $K_\fq$-algebra.] \,
\begin{itemize}
	\item The algebra is generated multiplicatively by six canonical basis elements: three
$T_i$ and three $D_i$, with $i$ defined modulo $3$. 
	\item The defining relations are
\begin{align}
	T_i\, T_{i+1}   &=  1 + \fq^{-1}\,\chi_1\, D_i + \fq^{-2}\, D_i^2 , &
	T_{i+1}\, T_i   &=  1 + \fq\,\chi_1\, D_i + \fq^{2}\, D_i^2 , \\
	D_i\, D_{i+1}   &=  1 + \fq^{-1}\, T_{i+1} , &
	D_{i+1}\, D_i   &=  1 + \fq\, T_{i+1} , \\
	T_i\, D_{i+1}   &=  \chi_1 + \fq^{-1}\, D_i + \fq\, D_{i-1} , &
	D_{i+1}\, T_i   &=  \chi_1 + \fq\, D_i + \fq^{-1}\, D_{i-1} ,
\end{align}
where $\chi_1$ is again the $SU(2)$ fundamental character, together with
\begin{equation}
	T_i\, D_i = \fq^{-2}\, D_i\, T_i , \qquad
	T_i\, D_{i-1} = \fq^{2}\, D_{i-1}\, T_i \, .
\end{equation}
	\item The remaining canonical basis elements are the monomials
	$\fq^{ab}\, T_i^{a}\, D_{i}^{b}$ or $\fq^{-ab}\, T_i^{a}\, D_{i-1}^{b}$, possibly dressed by the flavour characters $\chi_k$. 
	\item The automorphism $\rho(T_i) = T_{i+1}$, $\rho(D_i) = D_{i+1}$ fixing $\chi_k$ generates the $\bZ_3$ symmetry of $K_\fq([A_1,D_3])$.
	\item The trace is defined by the cyclicity axiom and the three elementary traces
	 \begin{align}\label{eq:a1d3trace}
	\Tr 1 &= \frac{1}{D}\sum_{n\ge0}(-1)^n\,\chi_{2n}(\mu)\,\fq^{\,3n(n+1)}
	= 1 + \chi_2\fq^2 + (\chi_0+\chi_2+\chi_4)\fq^4 + \ldots , \\
	\Tr D_0 &= \frac{1}{D}\sum_{p\ge1}(-1)^p\,\chi_{2p-1}(\mu)\,
	    \big(\fq^{\,p(3p-1)-1}-\fq^{\,p(3p+1)-1}\big) , \\
	\Tr T_0 &= \frac{1}{D}\sum_{m\ge0}(-1)^m\,\chi_{2m}(\mu)\,
	    \big(\fq^{\,3m(m+1)-1}-\fq^{\,m(3m+1)-1}-\fq^{\,(m+1)(3m+2)-1}\big) .
	\end{align}
	using the Weyl-Kac denominator
	\begin{equation}
	D(\mu,\fq) = (\fq^2;\fq^2)_\infty\,(\mu^2\fq^2;\fq^2)_\infty\,
	             (\mu^{-2}\fq^2;\fq^2)_\infty \, ,
	\end{equation}
\end{itemize}
\end{definition}
\begin{remark} The same considerations apply to this definition as for the case of $K_\fq([A_1,A_2])$. In particular, we can write recursion relations reducing all traces to the three elementary traces above. Negative powers of $\fq$ in the recursion then fix the trace uniquely up to the overall scale. \end{remark}

\begin{remark}
\label{rem:a1d3-miracle}
We also encounter an unexplained ``minor miracle''. Namely, the minors of the $2 \times 3$ matrix of linear constraints 
divided by the neat common factor $(\mu^2\fq^2;\fq^2)_\infty\,(\mu^{-2}\fq^2;\fq^2)_\infty$ reproduce directly
$\Tr 1$,$\Tr T_0$,$\Tr D_0$. 
\end{remark}

This $K_\fq$ algebra belongs to another infinite conjectural family $K_\fq([A_1,D_{2n+1}])$. More generally, 
there is an ADE family associated to cluster algebras with finite mutation orbits.
 
\section{K-theoretic Coulomb branch algebras as $K_\fq$ algebras}
\label{sec:coulomb}
The fusion algebras for supersymmetric gauge theories provide a large collection of $K_\fq$ algebras $K_\fq[G,M]$ labelled by a choice of gauge group $G$ and symplectic matter representation $M$. The definition of ``gauge group'' is a bit subtle in 4d gauge theory. It may take the form of an actual (compact, connected) global form $G$ for the underlying (reductive) gauge Lie algebra $\fg$, but there are other options we will describe below. A second subtlety is a restriction on possible $M$'s imposed by an anomaly cancellation condition. The anomaly always cancels when $M=T^*N$. In the rest of this Section we will assume this restriction on the matter content. It would be interesting to lift it, but some crucial definitions below depend on it. 

At the level of $K_\fq$ algebras, we can ignore the subtlety about global forms and define a ``minimal'' $K_\fq[\fg,M]$ algebra generated by line defects which are present for all choices of gauge group. The choice of global form allows one to (maximally) extend $K_\fq[\fg,M]$ to a larger $K_\fq[G,M]$ algebra which has $K_\fq[\fg,M]$ as a canonical sub-algebra. 

We begin our description of $K_\fq[G,M]$ as a $K_\fq$-algebra by labelling its canonical basis elements. Following  
\cite{Kapustin:2007wm,Kapustin:2006hi} we conjecture that the canonical basis elements $L_{m,e}$ in $K_\fq[\fg,M]$ will be labelled by a pair \begin{equation}
	(m,e) \;\in\; \big(\Lambda_m \times \Lambda_e\big)/\mathrm{Weyl} \,,
\end{equation} 
where $\Lambda_m$ is the magnetic (co)root lattice and $\Lambda_e$ is the root lattice, acted simultaneously by the Weyl group of $\fg$. If we need to pick canonical representatives, we can take $m$ to be dominant for $\fg$ and $e$ to be dominant for the Levi sub-algebra $\fg_m$ preserved by $m$.\footnote{The $L_{m,e}$ generators represent ``'t Hooft-Wilson lines'': $m$ represents the magnetic charge of a 't Hooft line defect, while $e$ is identified as the weight of an irreducible representation of $\fg_m$, selecting a Wilson line placed on top of the 't Hooft line defect.} 

A choice of global form $G$ of the gauge group enlarges the space of labels to 
\begin{equation}
	(m,e) \;\in\; \Lambda_G/\mathrm{Weyl} \,,
\end{equation} 
for some sub-lattice $\Lambda_G$ of $\Lambda^{w}_m \times \Lambda^{w}_e$, the product of magnetic (co)weights and weights lattices, Lagrangian under the natural pairing between $\Lambda^{w}_m/\Lambda_m$ and $\Lambda^{w}_e/\Lambda_e$ \cite{Gaiotto:2010be, Aharony:2013hda,Gaiotto:2014kfa}.  Conventional global forms correspond to lattices 
$\Lambda^G_m \times \Lambda^G_e$ of magnetic coweights and weights for $G$. 

Whenever the matter representation $N$ is non-trivial, we will be able to flavour $K_\fq[G,T^*N]$. We can always have a $U(1)$ favour factor for each irreducible summand of $N$. Multiple copies of the same irrep enhance $U(1)^n$ to $U(n)$. If the irrep is real or pseudo-real, the flavour group can be further enhanced to $USp(2n)$ or $SO(2n)$ respectively. If $N$ includes copies of both $R$ and $R^\vee$ there will be similar enhancements. Such further enhancement will not be manifest in the formulae we use in this Section, but can be typically restored with some extra work.  

The map $\rho$ acts on canonical labels as: 
\begin{equation}\label{eq:rho-witten}
	\rho : (m,e) \;\longmapsto\; (-m,\; \kappa(m)-e) \,,
\end{equation}
with $\kappa(m)$ determined by the difference of quadratic Casimirs of $\fg$ on $M$ and on the adjoint
representation:
\begin{equation}
\kappa(m) \;=\;
\sum_{\substack{\alpha \in \mathrm{Adj} \\ \langle \alpha,m\rangle > 0}}
   \langle \alpha,m\rangle\,\alpha
\;-\;
\sum_{\substack{w \,\in\, N \\ \langle w,m\rangle > 0}}
   \langle w,m\rangle\, w \, .
\end{equation}
In the presence of flavour, $\rho$ also acts on the the flavour charge $q_f$ by $q_f \to \kappa_f(m)-q_f$, with
\begin{equation}
\kappa_f(m) \;=\; \;-\;
\sum_{\substack{w \,\in\, N \\ \langle w,m\rangle > 0}}
   \langle w,m\rangle\, w_f\, ,
\end{equation}
where $w_f$ are the flavour weights of $N$. The shift is compatible with the $U(n)$ flavour enhancement for $n$ copies of the same irrep in $N$,
as it lives in the diagonal $U(1)$ Abelian subgroup. 
 
 \subsection{The Abelianized presentation.}
The algebra $K_\fq[G,T^*N]$ admits a faithful presentation as a sub-algebra in a certain rational variant of the quantum torus algebra \cite{Alday:2009fs,Drukker:2009id,Gomis:2009xg,Gomis:2010kv,Gomis:2011pf,Ito:2011ea,Gang:2012yr,Bullimore:2015lsa,Dedushenko:2018icp,Braverman:2016wma,Braverman:2016pwk,Finkelberg:2017qpy,Schrader:2019hzd,Cautis:2023jkq} associated to $\Lambda_G$. Concretely, we map algebra elements to finite collection of rational functions $f_m(v)$ of ``gauge fugacities'' $v$ valued in the Cartan torus on $G$ and of a magnetic coweight $m$, invariant under simultaneous Weyl transformations of both. More precisely, these functions take the form $v^e g_m(v^\alpha)$, with $\alpha$ denoting roots and $(m,e) \in \Lambda_G$. We define a bar involution on $f_m(v)$ in the obvious way, as $\fq \to \fq^{-1}$ with $v$ fixed. The image of $L_{m,e}$ will be bar-invariant.

The denominator factors of $f_m$ must take the form $(1-\fq^k v^\alpha)$ with integral $k$. The product is defined as 
\begin{equation}\label{eq:fgprod}
	(f \cdot g)_m(v) = \sum_{m'} f_{m'}\big(\fq^{\,m'-m} v\big)\;
	g_{m-m'}\big(\fq^{\,m'} v\big)\; \cocy_{m',\,m-m'}(v)
\end{equation}	
for a cocycle which includes one denominator factor per positive
root, and one numerator factor per $N$ weight:
\begin{align}\label{eq:ccprod}
	\cocy_{a,b}(v) \;=\;
	&\prod_{\alpha>0}\ \cocy^{\rm vec}
	  \big(\langle a,\alpha\rangle,\ \langle b,\alpha\rangle;\ v^{\alpha}\big) \cr
	\times\ &\prod_{i}\ \prod_{w\in N}
	  \cocy^{\rm hyp}\big(\langle a,w\rangle,\ \langle b,w\rangle;\ \mu^{w_f} v^{w}\big) \,,
\end{align}
the pure-gauge cocycle being the first product alone. 

A direction contributes only when the two charges pair with it in strictly
opposite signs; otherwise its factor is $1$. Write $A$, $B$ for the two pairings, and set
\begin{equation}
	p \;=\; \min\big(|A|,|B|\big) \,, \qquad r \;=\; \big||A|-|B|\big| \,,
\end{equation}
and $(x;\fq^{2})_p = \prod_{j=0}^{p-1}(1-\fq^{2j}x)$. Then, for $A>0>B$,
\begin{align}\label{eq:ccclosed}
	\cocy^{\rm vec}(A,B;z) &\;=\;
	\frac{\fq^{\,|AB|}\ z^{\,p}}
	     {\big(\fq^{\,r}z;\fq^{2}\big)_{p}\ \big(\fq^{\,r+2}z;\fq^{2}\big)_{p}}\,, \cr
	\cocy^{\rm hyp}(A,B;x) &\;=\; \big(-\fq^{\,r+1}x;\fq^{2}\big)_{p}\,,
\end{align}
and $A<0<B$ gives the bar image, $\fq\to\fq^{-1}$ with $v$ fixed.

An alternative compact notation introduces some $U_m$ symbols with $U_m v^e = \fq^{2 (m,e)} v^e U_m$ and appropriate product rule, 
\begin{equation}\label{eq:cocydef}
	U_a\, U_b \;=\; \cocy_{a,b}\big(\fq^{\,a+b}v\big)\; U_{a+b} \,.
\end{equation}
so we can write formally 
\begin{equation}
	f = \sum_m f_m(\fq^m v) U_m \, .
\end{equation}
A disadvantage of this compact notation is that it may include fractional powers of $\fq$ for certain $\Lambda_G$, which cancel out in 
the final results so the algebra is well-defined over $\bZ[\fq, \fq^{-1}]$. 

Both $\rho$, the trace and the pairing will be inherited from analogous structures on the rational quantum torus algebra:
\begin{equation}\label{eq:rhotorus}
	\big(\rho f\big)_{-m}(v,\mu) \;=\;
	v^{\kappa(m)}\ \mu^{\kappa_f(m)}\ f_m\big(v^{-1},\,\mu^{-1}\big) \,.
\end{equation}
and 
\begin{equation}\label{eq:measure}
	\Tr\, f \;=\; \frac{\big(\fq^2;\fq^2\big)_\infty^{2\,\mathrm{rk}\,\fg}}{|W|}
	\oint_{|v|=1} f_0(v,\mu)\
	\frac{\prod_{\alpha}\ \big(v^{\alpha};\fq^2\big)_\infty\,
	                      \big(\fq^{2}v^{\alpha};\fq^2\big)_\infty}
	     {\prod_{w\in N}\ \big(-\fq\mu^{w_f} v^{w};\fq^2\big)_\infty\,
	                      \big(-\fq\mu^{-w_f}v^{-w};\fq^2\big)_\infty} \,,
\end{equation}
and  
\begin{align}\label{eq:Iexplicit}
	I_{f,g} \;=\; &
	\frac{\big(\fq^2;\fq^2\big)_\infty^{2\,\mathrm{rk}\,\fg}}{|W|}
	\ \sum_{m}\ \oint_{|v|=1}\ f_m\big(v^{-1},\mu^{-1}\big)\ g_m(v,\mu) \cr
	\times\ &\frac{\prod_{\alpha}\
	   \big(\fq^{\,|\langle m,\alpha\rangle|}v^{\alpha};\fq^2\big)_\infty\,
	   \big(\fq^{\,2+|\langle m,\alpha\rangle|}v^{\alpha};\fq^2\big)_\infty}
	  {\prod_{w\in N}\
	   \big(-\fq^{\,1+|\langle m,w\rangle|}\mu^{w_f} v^{w};\fq^2\big)_\infty\,
	   \big(-\fq^{\,1+|\langle m,w\rangle|}\mu^{-w_f} v^{-w};\fq^2\big)_\infty} \,.
\end{align}
The expected relation 
\begin{equation}
 I_{f,g} = \Tr \rho(f) \,g = \Tr g \,\rho^{-1}(f) 
\end{equation}
holds at the level of integrands up to shifts $v \to \fq^{\pm m} v$ of the integration contour, thanks to cancellations between the measure and 
the $\cocy_{-m,m}$ cocycle. As a consequence, the identity could fail 
depending on the integrand pole structure. Conjecturally, the poles in $f$ and $g$ are always compatible with the required contour shifts.

In the following we will also employ a specific trivialization of the cocycle:
\begin{equation}\label{eq:cocha}
	\cocha_m(v) \;=\;
	(-\fq)^{\frac12\sum_{\langle\alpha,m\rangle>0}\langle\alpha,m\rangle}
	\prod_{\substack{\alpha \,:\, \langle\alpha,m\rangle>0}}\
	\prod_{\substack{k\ \mathrm{even}\\ -2\langle\alpha,m\rangle\,<\,k\,\le\,0}}
	\big(1-\fq^{\,k}v^{-\alpha}\big)^{-1}
	\ \prod_{w\in N}\
	\prod_{\substack{k\ \mathrm{odd}\\ 2\langle m,w\rangle\,<\,k\,<\,0}}
	\big(1+\fq^{\,k}\mu^{w_f} v^{w}\big)
\end{equation}
which allows an alternative presentation 
\begin{equation}
	f = \sum_m d_m(v) u^m \, ,
\end{equation}
with
\begin{equation}	
	d_a(v) = f_a(\fq^a v) \cocha_a(v) \, .
\end{equation}

The image of $L_{m,e}$ in the rational quantum torus algebra receives two contributions: the leading Weyl orbit
\begin{equation}
	f_{w\cdot m}(v) = w\cdot \chi^{\fg_m}_e(v) \, ,
\end{equation}
with $\chi^{\fg_m}_e(v)$ being the irreducible $\fg_m$ character of weight $e$, and a ``bubbling correction''.

We are ready to formulate our main conjecture: 
\begin{conjecture}[Coulomb branch canonical bases]\label{conj:abeKalgebra}
A consistent $K_\fq$-algebra $K_\fq[G,T^*N]$ can be defined uniquely by three conditions on the bubbling corrections:
\begin{enumerate}[label=(A\arabic*), ref=A\arabic*, leftmargin=*, itemsep=5pt]
  \item\label{ab:bar} The bubbling correction must be Weyl covariant and bar invariant.
  \item\label{ab:ofq} The bubbling correction must be $O(\fq)$.
  \item\label{ab:res} The whole $L_{m,e}$ must satisfy a residue cancellation rule: $d_a$ has simple poles only and 
    \begin{equation}\label{eq:star}
	\mathrm{Res}_{\,v^\alpha = \fq^{2l}}\ \big[\, d_a + d_b \,\big] \;=\; 0\,,
	\qquad b = s_\alpha(a) - l\,\alpha^\vee \,,
\end{equation}
	under any Weyl reflection $s_\alpha$.
\end{enumerate}
\end{conjecture}
Note: the first two axioms use the $f_a(v)$ data but the last uses the $d_a(v)$ data for the same candidate $L_{m,e}$. Also, the bubbling correction 
is supported in practice on Weyl orbits of magnetic weight subdominant to $m$, but this support condition appears to follow from the above axioms. Finally, the $d$'s have a pre-factor which may involve a fractional power of $-\fq$ in some cases. The same fractional power appears in front of all $d$'s in a given expression and thus does not affect the residue cancellation condition. 

\subsection{The $K_\fq$ algebra axioms}
Some aspects of our conjecture are straightforward, some appear rather non-trivial. The residue cancellation rule appeared before \cite{Klyuev:2026bic}
as a combinatorial characterization of the ``quantized K-theoretic Coulomb branch algebra'' defined by BFN for any pair $(G,N)$,
though it was mostly developed for the conventional Coulomb branch algebra. In general, the literature on K-theoretic Coulomb branch algebras is still somewhat under-developed. The residue cancellation rule can also be stated as the existence of a well-defined action of 
the rational quantum torus elements onto the space of $G$ characters $\chi_e(v)$ via $u^m \chi_e(v) = \chi_e(\fq^{2m} v)$.
We will review the physical origin of this requirement in Appendix \ref{app:physics}.  

It thus seems plausible that our conjecture could be simplified to the statement that the K-theoretic Coulomb branch algebra has a canonical basis of generators $L_{m,e}$ which are bar invariant and have $O(\fq)$ bubbling corrections. 

A second observation is that the $O(\fq)$ bubbling axiom together with the $I_{f,g}$ formula immediately implies the $\delta_{a,b} + O(\fq)$ 
$K_\fq$-algebra axiom, as the leading Weyl orbits already contributes $\delta_{a,b}$ and the bubbling corrections contribute $O(\fq)$. The 
\begin{equation}
 I_{f,g} = \Tr \rho(f) \,g = \Tr g \,\rho^{-1}(f) 
\end{equation}
relation is much trickier as it requires non-trivial pole cancellations. Fortunately, these cancellations are also known to occur in the K-theoretic Coulomb branch algebra, roughly thanks to the residue cancellation rule. 

Finally, we should recall that the quantized K-theoretic Coulomb branch algebra can be identified with the Grothendieck group of a mathematically well-defined category $\mathrm{KP}(G,N)$, which we expect to coincide with the abstract category $\mathrm{Lines}_{\bC^*}[G,N]$ of topological line defects in the holomorphic-topological twist of the super-symmetric gauge theory. The category admits structures lifting the product, the trace and $\rho$. The bar involution was only shown to exist when $N=0$, but is expected to exist for all $N$. The categorical analogue of the canonical basis has a natural definition as well. It thus may be possible to use that categorical presentation to prove directly that $K_\fq[G,M]$ satisfies the axioms of a $K_\fq$ algebra. It would be particularly interesting to understand the categorical origin of the $I_{a,b} = \delta_{a,b} + O(\fq)$ axiom.

\subsection{Pure $SU(2)$ and $SO(3)$.}
We can start from a neat rank $1$ example: $\fsl_2$. There are two conventional global forms, $SU(2)$ and $SO(3)$, 
and a third global form which is isomorphic to $SO(3)$ (but not canonically). We work in conventions where weights and coweights are
integers, so that $v^{\alpha}=v^{2}$, the pairing of the root with a coweight is
$\langle\alpha,m\rangle = m$, and the three lattices are
\begin{equation}\label{eq:su2so3lat}
	K_\fq[\fsl_2]:\ (m,e) \ \text{both even} \,,\qquad
	K_\fq[SU(2)]:\ m \ \text{even} \,,\qquad
	K_\fq[SO(3)]:\ e \ \text{even} \,.
\end{equation}
The third global form uses general $m$ and even $e-m$. The Dirac pairing is
$\tfrac12(me'-m'e)$.

The two ``minimal'' lines, which we cannot simultaneously employ, are the fundamental Wilson line $L_{0,1} = v + v^{-1}$ 
and the minuscule 't Hooft line $L_{1,0} = U_1 + U_{-1}$. The simplest lines shared by all theories are 
\begin{equation}\label{eq:su2thooft}
	L_{0,2} = v^{-2}+1+v^{2} \,,\qquad
	L_{2,0}= L_{1,0}^2 = U_2 + \frac{\big(\fq+\fq^{-1}\big)\,v^{2}}
	           {\big(1-\fq^{2}v^{2}\big)\big(1-\fq^{-2}v^{2}\big)} + U_{-2} \,.
\end{equation}
Here we cheated a bit: we squared $L_{1,0}$ and observed that the resulting bubbling correction is $O(\fq)$. But this example is simple enough that we can see the residue cancellation axiom at work. Passing to $d_a = f_a(\fq^a v)\cocha_a(v)$,
\begin{equation}\label{eq:su2dform}
	d_{-2} = \frac{-\fq}{\big(1-v^{2}\big)\big(1-\fq^{-2}v^{2}\big)} \,,\quad
	d_{0} = \frac{\big(\fq+\fq^{-1}\big)v^{2}}
	              {\big(1-\fq^{-2}v^{2}\big)\big(1-\fq^{2}v^{2}\big)} \,,\quad
	d_{2} = \frac{-\fq^{3}v^{4}}{\big(1-v^{2}\big)\big(1-\fq^{2}v^{2}\big)} \,.
\end{equation}
Three $v^{2} = \fq^{2l}$ residues matter, at $l=0,\pm1$, and each of the three $d$'s has a pole on
exactly two of them. The cancellation involves $a$ and $b=-a-2l$:
\begin{align}\label{eq:su2res}
	l=0,\ \ v^{2}=1:&\qquad d_{-2}+d_{2}
	  = \frac{-\fq + \big(\fq^{3}-\fq\big)v^{2} - \fq\,v^{4}}
	         {\big(1-\fq^{-2}v^{2}\big)\big(1-\fq^{2}v^{2}\big)} \,, \cr
	l=1,\ \ v^{2}=\fq^{2}:&\qquad d_{-2}+d_{0}
	  = \frac{-\fq + \big(\fq+\fq^{3}\big)v^{2}}
	         {\big(1-v^{2}\big)\big(1-\fq^{2}v^{2}\big)} \,, \cr
	l=-1,\ \ v^{2}=\fq^{-2}:&\qquad d_{0}+d_{2}
	  = \frac{\big(\fq+\fq^{-1}\big)v^{2}-\fq^{-1}v^{4}}
	         {\big(1-\fq^{-2}v^{2}\big)\big(1-v^{2}\big)} \,,
\end{align}
in each case the offending factor having disappeared from the denominator of the
sum.

Without the bubbling correction, the residue cancellation would obviously fail. The resulting $f_0$ is rigid. 
In particular, a Laurent polynomial in $v$ cannot be simultaneously bar invariant and $O(\fq)$. 
Similarly, we can e.g. compute
\begin{equation}\label{eq:su2dyonic}
	L_{2,1} \;=\; \fq\,v\,U_2
	 \;+\; \frac{v\big(1+v^{2}\big)}
	            {\big(1-\fq^{2}v^{2}\big)\big(1-\fq^{-2}v^{2}\big)}
	 \;+\; \fq\,v^{-1}U_{-2} \,,
\end{equation}
 
\subsection{$SU(2)$ and $SO(3)$ with adjoint matter.}
Taking $N$ to be the adjoint representation we obtain $N=2^*$ $SU(2)$ or $SO(3)$ gauge theories. 
Now $\kappa$ vanishes identically as the gauge and matter contributions cancel each other. Hence
\begin{equation}
	\rho : (m,e) \longmapsto (-m,-e) \,,\qquad \kappa_f(m) = -m \,.
\end{equation}

We still have the fundamental Wilson line $L_{0,1} = v + v^{-1}$ and the minuscule 't Hooft line $L_{1,0} = U_1 + U_{-1}$.
Squaring the latter and imposing the $O(\fq)$ axiom we find 
\begin{equation}\label{eq:adjsquare}
	L_{1,0}\cdot L_{1,0} \;=\; L_{2,0} + \mu \,,
\end{equation}
and 
\begin{equation}\label{eq:adjbubble}
	f_0(L_{2,0}) \;=\; \mu\ \frac{\big(\fq+\fq^{-1}\big)\big(\mu+\mu^{-1}\big)v^{2}
	                     + \big(1+v^{2}\big)^{2}}
	                    {\big(1-\fq^{2}v^{2}\big)\big(1-\fq^{-2}v^{2}\big)} \,.
\end{equation}

\section{RG flows of $K_\fq$ algebras}\label{sec:rg}
We now introduce the notion of RG flow $K_\fq[\UV] \to K_\fq[\IR]$ between $K_\fq$ algebras , axiomatizing the properties of certain ``Coulomb'' RG flows between the corresponding supersymmetric theories \cite{Ambrosino:2025qpy}. 
\begin{definition}[RG flow]
\label{def:rg}
The data of the RG flow consists of:
\begin{enumerate}[label=(R\arabic*), ref=R\arabic*, leftmargin=*, itemsep=5pt]
  \item\label{ax:lattice}
   	A lattice $\Gamma$ and a $K_\fq$ algebra $K_\fq[\IR]$ graded by $\Gamma$, e.g. canonical basis elements 
	carry a charge $\gamma \in \Gamma$, additive under multiplication. 
  \item\label{ax:map}
  A second $K_\fq$-algebra $K_\fq[\UV]$ equipped with an injective algebra map to $K_\fq[\IR]$:
	\begin{equation}
		\RG:\quad  L^\UV_a \to \RG_a \in K_\fq[\IR]\, ,
	\end{equation}
    where $\RG_a$ is bar-invariant and in particular an integral linear combination of $[n]_\fq L^\IR_b$ where $[n]_\fq$ are $\fq$-numbers.\footnote{A no exotic conjecture states that the coefficients are non-negative, but we will not impose it as an axiom.}
  \item\label{ax:S}
  A ``spectrum generator'' $S_\RG$, or just $S$, a formal element in $K_\fq[\IR]$ of the form
	\begin{equation}
		S \in \bZ[\fq,(1-\fq^{2n})^{-1}] \otimes K_\fq[\IR]_{\Gamma_+} 
	\end{equation}
where the support $\Gamma_+$ is a convex cone generated by a finite collection of simple charges.\footnote{This condition seems a bit artificial, but we have not encountered any example where it would need to be relaxed.} 
This satisfies an intertwining condition 
	\begin{equation}
		\RG_a S = S \rho^{-1}_\IR(\RG_{\rho_\UV(a)}) = L^\IR_{\ell(a)} + O(\fq)\, ,
	\end{equation}	
	where $\ell(a)$ is a bijective map from UV labels to IR labels. In particular, $S = 1+ 	O(\fq)$.
  \item\label{ax:rgtrace}
  Defining a notation $(L_a|L_b) = I_{a,b}$, we have
	\begin{equation}
		(L^\UV_a |L^\UV_b)_\UV = (\RG_a S | \RG_b S)_\IR \, .
	\end{equation}
\end{enumerate}
\end{definition}
\begin{remark} The $\RG_a S = L^\IR_{\ell(a)} + O(\fq)$ axiom fully characterizes $\RG_{\ell^{-1}(a)}$ given $S$ and $a$.  The full RG flow and UV $K_\fq$-algebra can thus be reconstructed from the IR $K_\fq$-algebra and $S$. \end{remark} 
This observation leads to an obvious question: which $S$ can be added to a given $K_\fq[\IR]$ to reconstruct a valid $K_\fq[\UV]$? The answer to this question appears to be amazingly rich, recovering and extending the theory of cluster algebras. 

\begin{remark} RG flows can be composed. Consider RG flows ``$1 \to 2$'' and ``$2 \to 3$'' from a $K_\fq$ algebra ``1'' to a $K_\fq$ algebra ``2'' and from a $K_\fq$ algebra ``2'' to a $K_\fq$ algebra ``3''.
We can easily obtain an RG flow ``$1 \to 3$'' from ``1'' to ``3''. The $\RG$ maps are composed in the usual way, 
\begin{equation}
	\RG^{1 \to 3} = \RG^{2\to 3} \circ \RG^{1 \to 2} \, ,
\end{equation}
and 
\begin{equation}
	S^{1 \to 3} = \RG^{2\to 3}(S^{1 \to 2}) \cdot S^{2 \to 3} \, .
\end{equation}
\end{remark}

\begin{remark} As $S$ is a formal element, it is potentially dangerous to re-write the inner products as traces. But $I_{a,b} = I_{\rho^{-1}(b),\rho^{-1}(a)}$ is safer and this operation and its inverse allows one to derive a family of new RG flows from any given one, by mapping  
\begin{align}
	(\RG,S) \to \left( \rho^{-1}_\IR(\RG_{\rho_\UV(a)}), \rho^{-1}_\IR(S) \right) \, .
\end{align}
\end{remark}
We can generalize this construction as follows. 

\begin{definition}[Wall-crossing]
\label{def:wallcrossing}
A wall-crossing relation between two RG flows $(\RG^{(1)},S^{1 \bar 1})$ and $(\RG^{(2)},S^{2 \bar 2})$ between the same pair of $K_\fq$-algebras is defined by a factorization  
\begin{equation}
	S^{1 \bar 1} = S^{12} S^{2 \bar 1} \qquad \qquad S^{2 \bar 2} = S^{2 \bar 1} \rho^{-1}_\IR(S^{12}) \,,
\end{equation}
such that
\begin{align}
	\RG^{(1)}_a S^{12} &= S^{12} \RG^{(2)}_a\cr
	\RG^{(2)}_a S^{2 \bar 1} &= S^{2 \bar 1}  \rho^{-1}_\IR(\RG^{(1)}_{\rho_\UV(a)}) \,.
\end{align}
\end{definition}
Note: the partial factors in $S$ do {\it not} necessarily have a good behaviour as $\fq \to 0$. 

\begin{remark} In a physical context, we have a linear central charge function $Z: \Gamma \to \bC$ and the positive cone $\Gamma_+$ satisfies 
$\mathrm{Im} \, Z_\gamma>0$. This condition can be rotated by a phase to $\mathrm{Im}\,  e^{- i \vartheta} Z_\gamma>0$
and we have a locally constant family of RG flows $(\RG^\vartheta,S^\vartheta)$ defined on the universal cover of $S^1$,
related by wall-crossing. The wall-crossing induced by $\rho^{\pm 1}$ corresponds to a shift by $\pm \pi$ of $\vartheta$. \end{remark} 

\begin{remark}
	There are natural RG flows from $K_\fq[G,N]$ to $K_\fq[G,\emptyset] \otimes Q_\fq[\Gamma]$, with $\Gamma = \bZ$:
	\begin{equation}
		S = \prod_{w \in N} E_\fq(\mu \,v^w)
	\end{equation} 
	re-interpreted as a linear combination of characters $\chi_w(v)$. Here $\mu$ is the generator of $Q_\fq[\Gamma]$. When combined with the conventional Abelianization formulae, this recovers a characterization of $K_\fq[G,N]$ by \cite{schrader2026ktheoreticcoulombbranchesquiver}. If $N$ has $n$ summands, we can introduce a separate flavour charge for each, $\Gamma = \bZ^n$ and 
	this RG flow factors in an obvious manner. 
\end{remark}

\subsection{Seiberg-Witten RG flows}
The best understood examples of RG flows of $K_\fq$-algebras land on the quantum torus: $K_\fq[\IR] = Q_\fq(\Gamma)$. 
In order to describe them in our language, we need to introduce some extra notation. 

Our general RG flow axioms require $S$ to belong to the 
group $\cG$ of elements of the form
\begin{equation}
	1+\sum_{\gamma\in\Gamma_+} g_\gamma(\fq)X_\gamma \qquad \qquad g_\gamma(\fq)\in \bZ\left[\fq,\fq^{-1},(1-\fq^{2n})^{-1}|_{n>0}\right] .
\end{equation}
This can be formalized as a pro-nilpotent group. For example, our assumption that $\Gamma_+$ is generated by a collection of simple charges $\gamma_i$ means that the coefficients in the product of two such elements are finite sums. 

We introduce the quantum dilogarithm 
\begin{equation}
	E_\fq(x)=(-\fq x;\fq^2)^{-1}_\infty = \sum_{n=0}^\infty \frac{(-\fq)^n}{(\fq^2)_n} x^n = 1 - \frac{\fq}{1-\fq^2} x + \cdots
\end{equation}
and a generalization for $s\in \tfrac12\mathbb{N}$:
\begin{equation}\label{eq:mult}
  E^{(s)}_\fq(x)=\prod_{j=-s}^{s}
    E_\fq\bigl((-1)^{2s}\fq^{2j}x\bigr)^{(-1)^{2s}} \,.
\end{equation}
\begin{definition}We define the subgroup $\cE$ of $\cG$ generated by $E^{(s)}_\fq(X_\gamma)^{\pm 1}$ for $\gamma \in \Gamma_+$.  
\end{definition}

\begin{remark} It is easy to argue that every element $g \in \cE$ satisfies $g^{-1}(\fq) = g(\fq^{-1})$, but $\cE$ is smaller than the group 
of elements which satisfy that constraint. For example, it does not include $E_{\fq^2}(x)$. \end{remark}

\begin{definition}[Seiberg-Witten RG flow]
\label{def:sw}
We define a Seiberg-Witten RG flow as an RG flow with $K_\fq[\IR] = Q_\fq(\Gamma)$ and $S \in \cE$. A ``BPS quiver RG flow'' additionally satisfies:
\begin{equation} \label{eq:quiver}
	S = 1 - \fq \sum_i X_{\gamma_i} + O(\fq^2) \, .
\end{equation}
\end{definition}
We then propose:
\begin{conjecture}[Uniqueness of $S$]\label{conj:S-unique}
	There is an unique element $S \in \cE$ which satisfies the BPS quiver constraint (\ref{eq:quiver}).
\end{conjecture}
We believe this conjecture follows from a remarkable combinatorial property of $\cE$:
\begin{conjecture}[PBW factorization]\label{conj:pbw}
	Every element $g\in \cE$ can be factorized as 
	\begin{equation}
		g = \prod_{\gamma,s} E^{(s)}_\fq(X_\gamma)^{\Omega(\gamma,s)}
	\end{equation}
	where the product is taken in any chosen total order on 
	$\Gamma_+ \times \tfrac12\mathbb{N}$ and the exponents are integral.
\end{conjecture}
\begin{remark} If the factorization exists, is unique: we can proceed along the filtration of $\Gamma_+$ and at each step 
$\Omega(\gamma,s)$ is uniquely determined. Furthermore, it is easy to see that for every choice of total order 
we can always build a solution $S$ of equation (\ref{eq:quiver}) in that order, starting with $\prod_i E_\fq(X_{\gamma_i})$ and 
moving along the filtration: at each step we can use $\Omega(\gamma,s)$ to cancel the non-positive powers of $\fq$ in the coefficient of $X_\gamma$. If the factorization conjecture holds, $S$ is independent of the chosen total order. \end{remark}

\begin{remark} In a physical context, one only encounters the total order determined by $\arg Z_\gamma$. More general choices of total order seem mathematically consistent, though. \end{remark}

The ``BPS quiver'' in the moniker is a quiver $Q$ with nodes labelled by the $\gamma_i$, no self-loops and no two-cycles
arrows determined by 
\begin{equation}
	\langle \gamma_i,\gamma_j\rangle = \#(i \to j) - \#(j \to i) \, .
\end{equation}
The spectrum generator $S$ as an element of the sub-lattice generated by the node charges $\gamma_i$  only depends on $Q$, 
so we could denote it as $S[Q]$. It should be compared with other analogous objects associated to BPS quivers by the theory of cluster algebras or 
Hall algebras:
\begin{enumerate}
	\item For many BPS quivers it is possible to find a ``maximal green
	sequence'' of mutations, which can be converted into a finite sequence
	of charges $\gamma(a)$ in $\Gamma_+$ and then into a candidate
	spectrum generator:
		\begin{equation}
			S_{\mathrm{cluster}} \equiv \prod_a E_\fq(X_{\gamma(a)}) \, ,
		\end{equation}
		the product being ordered by increasing $a$, i.e.\ in the order of the sequence.
	Although clearly $S_{\mathrm{cluster}} \in \cE$, it is far from obvious that it satisfies (\ref{eq:quiver}).
	In all examples we considered, though, this is the case. We thus conjecture
	\begin{conjecture}[$S_{\mathrm{cluster}}=S$]\label{conj:cluster}
	\begin{equation}	
		S_{\mathrm{cluster}} = S \, ,
	\end{equation}
	when a maximal green sequence can be found. 
	\end{conjecture}
	\item In the special case of BPS quivers associated to character varieties, one can define an 
	$S_{\mathrm{DT}}$ element geometrically. 
	\begin{conjecture}[$S_{\mathrm{DT}}=S$]\label{conj:dt}
	\begin{equation}	
		S_{\mathrm{DT}} = S \, ,
	\end{equation}
	when the former can be defined. 
	\end{conjecture}
	It is far from obvious that $S_{\mathrm{DT}}$ satisfies (\ref{eq:quiver}).
	\item If the quiver is equipped with a potential $W$, one can define a ``Cohomological Hall Algebra'' $H(Q,W)$ whose equivariant character 
	defines another special element $S^{-1}_{\mathrm{CoHA}} = \chi(H(Q,W))$. This element depends on the choice of $W$. 
	\begin{conjecture}[$S_{\mathrm{CoHA}}=S$]\label{conj:coha}
	For any BPS quiver $Q$ there exist a choice of (infinitely mutable) $W$ such that 
	\begin{equation}	
		S_{\mathrm{CoHA}}= S \, ,
	\end{equation}
	\end{conjecture}
	Again, it is far from obvious that $S_{\mathrm{CoHA}}$ satisfies (\ref{eq:quiver}).
\end{enumerate}

Given $S(Q)$, we can solve the equations $\RG_\gamma S = X_\gamma + O(\fq)$ for all $\gamma$. The solution $\RG_\gamma$ exists and is unique as a bar-invariant formal sum supported on $\gamma + \Gamma_+$ with coefficients in $\bZ[\fq, \fq^{-1}]$. 
\begin{conjecture}[Upper cluster algebra]\label{conj:upper}
	For all BPS quivers $Q$, possibly embedded in a lattice $\Gamma$ which is bigger than the span of $\gamma_i$, 
	the $\RG_\gamma$ are actual elements in $Q_\fq(\Gamma)$ and define a sub-algebra over $\bZ[\fq, \fq^{-1}]$. There is a label permutation 
	$\rho_\UV$ which makes this data into a $K_\fq$ algebra $K_\fq(Q,\Gamma)$ with canonical RG flow to $Q_\fq(\Gamma)$.
	As an algebra, $K_\fq(Q,\Gamma)$ coincides with the ``upper cluster algebra'' of mutation-covariant Laurent polynomials in $Q_\fq(\Gamma)$ 	attached to $(\Gamma,Q)$ by the theory of cluster algebras.
\end{conjecture}
We tested this conjecture extensively. 

Finally, mutations of a BPS quiver should lead to wall-crossing transformations. More precisely, 
\begin{conjecture}[$S$ Mutation]\label{conj:mut} Let $Q$ be a BPS quiver with node charges $\gamma_i$, and $Q'=\mu_j(Q)$ the result of a mutation at the node $j$, aka the BPS quiver 
with node charges
\begin{equation}\label{eq:mut-charges}
	\gamma'_j = -\gamma_j , \qquad
	\gamma'_i = \gamma_i + [\langle \gamma_i, \gamma_j \rangle]_+ \, \gamma_j \quad (i \neq j) \, ,
\end{equation}
where $[x]_+ = \max(x,0)$.
Then the associated spectrum generators satisfy a wall-crossing relation of the form:
\begin{align}
	S &= E_\fq(X_{\gamma_j}) S'' \cr
	S' &= S'' E_\fq(X_{-\gamma_j}) \, ,
\end{align}
for an element $S''$ in the intersection of the associated groups $\cE$ and $\cE'$.
\end{conjecture}

\subsection{Factored RG flows}
Consider now a pair $(Q,Q')$ of a BPS quiver $Q$ and a full sub-quiver $Q'$. Embed $Q$ in a lattice $\Gamma$ for convenience. By picking an appropriate total order, we can factor
\begin{equation}
	S_Q = R_Q^{Q'} S_{Q'} 
\end{equation}
\begin{conjecture}[Factored RG flows]\label{conj:factored}
	We conjecture that $R_Q^{Q'}$ is in the image of $\RG^{Q',\Gamma}$:
	\begin{equation}
	S_Q = \RG^{Q',\Gamma}(S_Q^{Q'}) S_{Q'} \, ,
	\end{equation}
	for a formal element $S_Q^{Q'}$ in $K_\fq(Q',\Gamma)$ and $\RG^{Q,\Gamma}$ factors as $\RG^{Q',\Gamma}\circ \RG^{Q,Q',\Gamma}$ as well, defining an RG flow of $K_\fq$ algebras from $K_\fq(Q,\Gamma)$ to $K_\fq(Q',\Gamma)$.
\end{conjecture}

This construction is our main source of examples of RG flows between known $K_\fq$ algebras. It would be interesting to 
characterize which pronilpotent group of formal sums in $K_\fq(Q',\Gamma)$ is the natural place for such $S_Q^{Q'}$ RG flow generating elements.

\subsection{Flavoured BPS quivers}
When a BPS quiver has $N$ identical nodes, the associated $K_\fq$ algebra is expected to be $SU(N)$-flavoured. 
This can be made explicit by a slight modification of our axioms. 
Introduce a flavour group $G_f$ and a flavoured quantum dilogarithm for a finite-dimensional irrep $R$ of $G_f$:
\begin{equation}\label{eq:fleq}
  E^{(s;R)}_\fq(x)=\prod_{w \in R} E^{(s)}_\fq(\mu^w x) \,.
\end{equation}
expanded into characters $\chi_w(\mu)$ of $G_f$ irreps. 
\begin{definition}We define the group $\cE_{G_f}$ generated by $E^{(s;R)}_\fq(X_\gamma)^{\pm 1}$ for $\gamma \in \Gamma_+$.  
\end{definition}
We can extend the PBW factorization conjecture \ref{conj:pbw} to $\cE_{G_f}$ in the obvious way. 

The most general situation has every $\gamma_i$ in the BPS quiver appear with multiplicity $N_i$, leading to 
$G_f = \prod_i SU(N_i)$ and an axiom 
\begin{equation} \label{eq:flavouredquiver}
	S = 1 - \fq \sum_i \chi^{\square_i}(\mu) X_{\gamma_i} + O(\fq^2) \, ,
\end{equation}
where $\chi^{\square_i}(\mu)$ are the characters for the fundamental representation of $SU(N_i)$. 

We can then define an RG flow of flavoured KAlgebras landing on the combination of the representation ring of $G_f$ and 
$Q_\fq(\Gamma)$. This is a natural source of flavoured $K_\fq$-algebras. 

\section{Skein algebras as $K_\fq$-algebras}\label{sec:skein}
This section explores briefly the algebraic consequences of the existence of 4d ${\cal N}=2$ supersymmetric theories of ``class S'' \cite{Witten:1997sc,Gaiotto:2009we,Gaiotto:2009hg}
and associated $K_\fq$ algebras. An important source of insights is to look at structures associated to the underlying 6d $(2,0)$ superconformal field theories. In particular, the line defects we consider in 4d theories of class S arise from supersymmetric, half-BPS 2d defects in 6d.

Our general expectation is that the physical properties of such 2d defects have implications on the structure of the category $\mathrm{Rep}_\fq(\fg)$ 
of finite-dimensional representations of $U_\fq(\fg)$, at least for simply-laced $\fg$. Intuitively, each irreducible representation $V_\lambda$ is associated to one such 2d defect $\mathbb{D}_\lambda$.   

Given a cyclic collection $(\lambda_a)$ of such defects, we can consider a 1d {\it junction} where the 2d defects collectively end. Such a junction can also be supersymmetric: the 6d theory has $16$ supercharges, the defect preserve $8$ and a junction can preserve $4$. This is the same amount of supersymmetry preserved by the line defects we study implicitly in the rest of the paper. Just as in 4d, the space transverse to the 
defects and junction is $\bR^3$. We propose to identify the space of $U_\fq(G)$ invariants
\begin{equation}
	\mathrm{In}_{(\lambda_a)} \equiv \mathrm{Inv}\left(V_{\lambda_1} \otimes \cdots \otimes V_{\lambda_n} \right)
\end{equation}
as a vector space over $\bZ[\fq,\fq^{-1}]$ with the Grothendieck group of such junctions, with the physical supersymmetric junctions 
defining a bar-invariant canonical basis. Accordingly, we conjecture
\begin{conjecture}[Canonical junctions]\label{conj:junctions}
The spaces $\mathrm{In}_{(\lambda_a)}$ have bar-invariant canonical bases over $\bZ[\fq,\fq^{-1}]$ such that 
\begin{enumerate}
	\item Cyclic rotation of the $\lambda_a$ acts as a permutation.
	\item The collapse of any webs of junctions to a single junction defines linear operations over $\bZ[\fq,\fq^{-1}]$ with coefficients which are palindromic, e.g. invariant under $\fq \to \fq^{-1}$. 
\end{enumerate}
\end{conjecture}
\begin{remark} The additional naive conjecture that the coefficients would be non-negative linear combinations of $[n]_\fq$ numbers appears to be false
in examples. \end{remark}
\begin{remark} This conjecture is proven in upcoming work \cite{BK} announced to the author in a private communication.\end{remark}
Given such a conjectural set of canonical junctions, we can now propose a $K_\fq$ algebra structure on {\it skein algebras} $\mathrm{Sk}_\fq(\fg,C)$ associated to a Riemann surface $C$. Our proposal should admit extensions to Riemann surfaces with appropriate decorations \cite{Jordan:2021hop}. Recall that the skein algebra is the operator algebra of Chern-Simons theory compactified on $C$. Concretely, generators are skeins of  $\mathrm{Rep}_\fq(\fg)$ representations, either forming closed curves or more complicated webs with invariant junctions, drawn in $C \times [-1,1]$. In order for the skein algebra to be defined over $\bZ[\fq,\fq^{-1}]$ we need to impose some restrictions on skein labels. If we restrict to labels in the root lattice 
we will obtain a ``small'' skein algebra analogous to $K_\fq[\fg,N]$ in the gauge theory story. Weaker restrictions labelled by a Lagrangian submanifold of $H^1(C,Z(\fg))$ are possible, analogous to general global forms $K_\fq[G,N]$ in the gauge theory story.

The skein algebra is equipped with a bar involution defined by a reflection of the $[-1,1]$ segment. We now define a conjectural canonical basis in 
$\mathrm{Sk}_\fq(\fg,C)$: irreducible skeins drawn at $C \times \{0\}$, with the canonical junctions defined above. Skeins which are not irreducible 
will reduce to linear combinations of the irreducible skeins with $\fq$-number coefficients. We define $\rho$ as a dualization operation on all skeins,
assuming some action can be defined on the canonical junctions. 

This leaves us with the problem of defining a trace. This is more challenging. Assume for simplicity that $\rho^2=1$, which is the case in the absence of irregular punctures, so we have a conventional trace. An universal source of traces is the dual of the skein module of $S^1 \times C$, which describe traces which can be successfully extended to functions on all skeins drawn on 
$S^1 \times C$. In turn, a way to produce interesting elements of that vector space is to specify a four-manifold bounded by $S^1 \times C$: the path integral of Kapustin-Witten theory on such a manifold with Neumann boundary conditions produces in particular a trace on $C$. We conjecture that the trace we are after is associated to $D^2 \times C$. 

\section*{Acknowledgements}
We would like to thank K. Costello, G. Moore, B. Webster, J. Kamnitzer, H. Williams, Y. Soibelman for useful conversations and feedback on the draft.
This research was supported in part by a grant from the Krembil Foundation. DG is supported by
the NSERC Discovery Grant program and by the Perimeter Institute for Theoretical Physics.
Research at Perimeter Institute is supported in part by the Government of Canada through the Department of Innovation, Science and Economic
Development Canada and by the Province of Ontario through the Ministry of Colleges and Universities.

\appendix

\section{Physical motivations}\label{app:physics}
We refer to \cite{Gaiotto:2010be, Kapustin:2007wm,Cordova:2016uwk,Ambrosino:2025qpy,Gaiotto:2024fso} for explanations of the physical origin of most of the algebraic structures in this paper. Our focus are the constraints on the $\fq \to 0$ behaviour we propose of several quantities of interest in the bulk of the paper. 

All of these statements concern the behaviour of ``Schur indices'', which can be defined as protected equivariant characters for the cohomology of an appropriate nilpotent ``holomorphic-topological'' supercharge acting on local operators \cite{Kapustin:2006hi,Niu:2021jet}. The original definition of the Schur index was as a specialization of the super-conformal index \cite{Kinney:2005ej,Romelsberger:2005eg,Gadde:2011ik,Gadde:2011uv} and thus only applied to super-conformal SQFTs, aka SCFTs. The general super-conformal index can also be interpreted as the protected equivariant character for the cohomology of an ``holomorphic'' super-charge, but the definition requires the full $U(1)_r \times SU(2)_R$ R-symmetry of the theory to be non-anomalous, and thus cannot be extended to non-conformal 4d ${\cal N}=2$ SQFTs: their $U(1)_r$ is always anomalous. 
The Schur specialization only requires $SU(2)_R$, which is a symmetry of all known  SQFTs with vanishing FI parameters. 

The conventional Schur index computes $I_{1,1}$ in our language. The condition $I_{1,1} = 1 + O(\fq)$ is easily understood for SQFTs: the super-conformal index is a graded trace over 
the space of states on a sphere, which is an unitary representation of the super-conformal group $SU(2,2|2)$. The power of $\fq$ is non-negative because of a straightforward unitarity bound,
which is only saturated by the identity operator. Hence the $\fq\to 0$ condition. In the conformal case, the Schur index can be generalized to count local operators 
on a super-conformal line defect, which breaks the super-conformal group to an $OSp(4^*|2)$. The same unitarity considerations apply and imply $I_{a,b} = \delta_{a,b}+ O(\fq)$
for super-conformal $L_a$ and $L_b$. Concretely, the super-algebra is
\begin{align}
[D,P] &= P, \qquad [D,K] = -K, \qquad [K,P] = 2D, \\[2pt]
[M_{\alpha\beta},M_{\gamma\delta}] &= \tfrac12\big(
  \epsilon_{\gamma\alpha}M_{\beta\delta}+\epsilon_{\gamma\beta}M_{\alpha\delta}
 +\epsilon_{\delta\alpha}M_{\beta\gamma}+\epsilon_{\delta\beta}M_{\alpha\gamma}\big), \\[2pt]
[R_{IJ},R_{KL}] &= \tfrac12\big(
  \epsilon_{KI}R_{JL}+\epsilon_{KJ}R_{IL}
 +\epsilon_{LI}R_{JK}+\epsilon_{LJ}R_{IK}\big), \\[4pt]
[D,\mathcal{Q}_{\alpha I}] &= \tfrac12\,\mathcal{Q}_{\alpha I}, \qquad
[D,\mathcal{S}_{\alpha I}] = -\tfrac12\,\mathcal{S}_{\alpha I}, \\[2pt]
[K,\mathcal{Q}_{\alpha I}] &= \mathcal{S}_{\alpha I}, \qquad\ \
[P,\mathcal{S}_{\alpha I}] = -\mathcal{Q}_{\alpha I}, \qquad
[P,\mathcal{Q}] = [K,\mathcal{S}] = 0, \\[2pt]
[M_{\alpha\beta},\mathcal{Q}_{\gamma I}] &= \tfrac12\big(
  \epsilon_{\gamma\alpha}\mathcal{Q}_{\beta I}+\epsilon_{\gamma\beta}\mathcal{Q}_{\alpha I}\big),
\qquad
[R_{IJ},\mathcal{Q}_{\alpha K}] = \tfrac12\big(
  \epsilon_{KI}\mathcal{Q}_{\alpha J}+\epsilon_{KJ}\mathcal{Q}_{\alpha I}\big), \\[4pt]
\{\mathcal{Q}_{\alpha I},\mathcal{Q}_{\beta J}\} &= \epsilon_{\alpha\beta}\epsilon_{IJ}\,P,
\qquad
\{\mathcal{S}_{\alpha I},\mathcal{S}_{\beta J}\} = \epsilon_{\alpha\beta}\epsilon_{IJ}\,K, \\[2pt]
\{\mathcal{Q}_{\alpha I},\mathcal{S}_{\beta J}\} &=
  \epsilon_{\alpha\beta}\epsilon_{IJ}\,D
+ \epsilon_{IJ}\,M_{\alpha\beta}
- 2\,\epsilon_{\alpha\beta}\,R_{IJ}.
\end{align}
Here $M$ are rotations, $R$ the R-symmetry generators, $P,D,K$ the 1d conformal generators, $Q$ the supercharges and $S$ the super-conformal partners. The Schur index is 
graded by something like $\fq^{2D - 2R_0}$, and the BPS bound gives $D \geq j + 2 r$ for operators of spin $j$ and $R$-symmetry spin $r$. 

We have not been able to identify a clear physical justification to extend such claims to non-conformal theories and line defects. For completeness, we will discuss some promising, but incomplete arguments. One possible route is to realize the Schur index $I_{a,b}$ of 
a non-conformal system as a graded trace over the space of states on a sphere with the help of some rigid supergravity background as in \cite{Festuccia:2011ws,Closset:2013vra}. Such a background will only preserve a ``super-isometry'' subgroup of $OSp(4^*|2)$, extending the isometry group of $S^3 \times S^1$ preserving the poles of $S^3$, aka the location of the line defects. For example, that could be $U(2|1)$. The presence of a decoupled $U(1)$ factor prevents any direct unitarity argument, though. At best we can write the exponent of $\fq$ as a linear combination of a generator bound by unitarity bounds and a generator which is quantized. It may be possible to further deform the system to a situation where the quantized generator can be bound, and transport this information back using the fact that continuously varying integers are necessarily constant. 

Other conjectures in the paper constrain the $\fq \to 0$ behaviour of other quantities which can be defined as Schur indices in the presence of interfaces and line defects. For example, the inner product 
\begin{equation}
		(L^{\IR}_a,RG_b S)_{\IR}
\end{equation}
counts local operators at a ``RG interface'' \cite{Cordova:2016uwk,Dimofte:2013lba} between the IR and UV theories, at a point where 
$L^{\IR}_a$ and $L^{\UV}_b$ simultaneously end. We also expect 
\begin{equation}\label{eq:int}
	\big(\fq^2;\fq^2\big)_\infty^{2\,\mathrm{rk}\,\fg}
	\oint_{|v|=1}\
	\frac{ v^{-e} \ f_m(v,\mu)[L_a] \prod_{\alpha}\
	   \big(\fq^{\,2+|\langle m,\alpha\rangle|}v^{\alpha};\fq^2\big)_\infty}
	  {\prod_{w\in N}\
	   \big(-\fq^{\,1+|\langle m,w\rangle|}\mu^{w_f} v^{w};\fq^2\big)_\infty\,
	   \big(-\fq^{\,1+|\langle m,w\rangle|}\mu^{-w_f} v^{-w};\fq^2\big)_\infty} \,.
\end{equation}
to count local operators at an interface between a $(G,N)$ and $(T_G,N)$ gauge theories, trivial for the matter fields and with the gauge group reduced to the diagonal Cartan $T_G$ at the interface, at a point where $L^G_a$ and $L^{T_G}_{m,e}$ also end. If we could argue that these expressions must only contain non-negative powers of $\fq$, we could derive many of the $O(\fq)$ claims in the man text. Unfortunately, even when conformal symmetry is present such configurations of intersecting defects can preserve at most $U(2|1)$ and we lack a unitarity constraint on the overall $U(1)$ factor. 

\section{$K_\fq$-algebras of finite type}\label{app:finite}
In this Appendix we review a collection of $K_\fq$-algebras which refine the notion of 
cluster algebras of finite type. These algebras have an ADE classification 
and are expected to be associated to certain theories of Argyres-Douglas type, 
denoted as $[A_1,\mathrm{ADE}]$. In practice, the finite type condition means that one can 
define the algebras in terms of a finite collection of multiplicative generators $L_a$, 
some of which ``$\fq$-commute'' pairwise:
\begin{equation}
	L_a L_b = \fq^{2(a,b)} L_b L_a \,,
\end{equation}
and can be assembled into the more complicated canonical basis elements
\begin{equation}
	L_{(a_1,\cdots, a_n)} \equiv \fq^{-\sum_{i<j} (a_i,a_j)} L_{a_1}\cdots L_{a_n} \,.
\end{equation}
All canonical basis elements are of this form, and are all distinct. The whole algebra can be presented by specifying the remaining $L_a L_b$ products for 
pairs which do not $\fq$-commute. 

The most ``economical'' presentation of the $K_\fq([A_1,ADE])$ algebras is of course to set up the cluster algebra machinery or its $Q_\fq$-charts analogue,
allowing the whole algebra to unfold from the combinatorics of cluster mutations. A more explicit presentation, though, is useful in order to 
assess the implications of the $I_{a,b}$ axioms, verify that they determine $I_{a,b}$ uniquely up to $1+O(\fq)$ factors and that they are self-consistent in 
these examples. 

\begin{table}[h]
\centering
\begin{tabular}{c|c|c}
\hline
$[A_1, ADE]$ & flavour symmetry & order of $\rho$ \\
\hline
$[A_1, A_{2n}]$   & ---    & $2n+3$ \\
$[A_1, A_{3}]$ & $SU(2)$ & $3$ \\
$[A_1, A_{2n+3}]$ & $U(1)$ & $2n+6$ \\
\hline
$[A_1, D_4]$              & $SU(3)$            & $4$ \\
$[A_1, D_{2n+1}]$         & $SU(2)$            & $2n+1$ \\
$[A_1, D_{2n+2}]$  & $SU(2)\times U(1)$ & $2n+2$ \\
\hline
$[A_1, E_6]$ & ---    & $14$ \\
$[A_1, E_7]$ & $U(1)$ & $10$ \\
$[A_1, E_8]$ & ---    & $16$ \\
\hline
\end{tabular}
\caption{Flavour symmetry and the order of the half-monodromy $\rho$ for the
$K_\fq([A_1, ADE])$ algebras of finite type. Here $n$ is a positive integer.}
\label{tab:finite_type}
\end{table}

The $[A_1,A_k]$ and $[A_1,D_k]$ AD theories have a class $S$ realization of type $A_1$, which allows a
presentations of the corresponding $K_\fq$ algebras as certain $A_1$ skein algebras. It is also possible to identify $[A_1,E_6]\simeq [A_2,A_3]$ and 
$[A_1,E_8]\simeq [A_2,A_4]$ as class $S$ theories of type $A_2$. 

\subsection{The odd polygon algebras, aka $K_\fq([A_1,A_{2k}])$.}\label{app:a1a2k}
The algebra $K_\fq([A_1,A_{2k}])$ is multiplicatively generated by
$L_{a;i}$ associated to the diagonals of a $(2k+3)$-gon. Here $(a;i)$
label the diagonal between the $i$-th and $i+a+1$-th vertices. The
edges of the polygon are all formally associated to the identity $1$.
The automorphism $\rho$ acts as $L_{a;i} \to L_{a;i+1}$ modulo $2k+3$
and generates a $\bZ_{2k+3}$ discrete symmetry.

It will be convenient to also write $L_{\overline{pq}} := L_{q-p-1;\,p}$
for the chord with endpoints $p, q \in \bZ/(2k{+}3)$, where $q-p$ is
taken in $\{2, \ldots, k+1\}$ modulo $2k{+}3$. The subscript
$\overline{pq}$ is the unordered pair $\{p, q\}$, so
$L_{\overline{pq}} = L_{\overline{qp}}$; edges
($q - p \equiv \pm 1$) are identified with $1$.

\begin{figure}[h]
  \centering
  \begin{tikzpicture}[scale=2]
    \foreach \i in {0,...,6} {
      \coordinate (V\i) at ({90 - 360*\i/7}:1);
      \fill (V\i) circle (1.5pt);
      \node at ({90 - 360*\i/7}:1.18) {$\i$};
    }
    \draw[gray!40] (V0) -- (V1) -- (V2) -- (V3) -- (V4) -- (V5) -- (V6) -- cycle;
    \draw[thick] (V0) -- (V1) node[midway,above,sloped] {\tiny $1$};
    \draw[thick] (V0) -- (V6) node[midway,above,sloped] {\tiny $1$};
    \draw[thick, blue] (V0) -- (V2) node[midway,above,sloped] {\tiny $L_{1;0}$};
    \draw[thick, blue] (V0) -- (V5) node[midway,above,sloped] {\tiny $L_{1;5}$};
    \draw[thick, red] (V0) -- (V3) node[midway,above,sloped] {\tiny $L_{2;0}$};
    \draw[thick, red] (V0) -- (V4) node[midway,above,sloped] {\tiny $L_{2;4}$};
  \end{tikzpicture}
  \caption{All chords incident to vertex $0$ on the heptagon ($k=2$):
    two length-$1$ edges (black, identified with the identity $1$),
    two length-$2$ short chords $L_{1;0}$ and $L_{1;5}$ (blue,
    orbit $a=1$), and two length-$3$ long diagonals $L_{2;0}$ and
    $L_{2;4}$ (red, orbit $a=2$). A chord is unoriented, so
    $L_{1;5} = L_{\overline{50}}$ has endpoints $(5, 5+2) = (5, 0)$
    and $L_{2;4} = L_{\overline{40}}$ has endpoints $(4, 4+3) = (4, 0)$.}
  \label{fig:heptagon-orbits}
\end{figure}
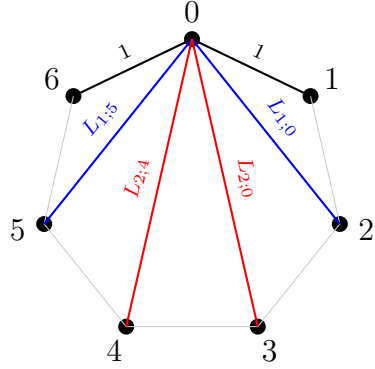

These generators satisfy ``quantum Ptolemy'' quadratic relations,
governed entirely by the geometry of the two chords involved.

\paragraph{Geometric chord pair.} Let $L_a, L_b$ be two letters with
endpoints $\{p_a, q_a\}$ and $\{p_b, q_b\}$. They are
\begin{itemize}
  \item[(D)] \emph{disjoint}: $\{p_a, q_a\} \cap \{p_b, q_b\} = \emptyset$
    and the chords do not cross,
  \item[(S)] \emph{vertex-sharing}: share exactly one endpoint,
  \item[(X)] \emph{crossing}: interiors meet transversally.
\end{itemize}

\paragraph{(D) and (S): $q$-commutation.}
For non-crossing pairs the product is a single canonical-basis monomial
and the two factors $q$-commute,
\begin{equation}
  L_a \cdot L_b \;=\; \fq^{2c_{ab}}\, L_b \cdot L_a,
  \qquad c_{ab} \in \{0, \pm 1\}, \quad c_{ba} = -c_{ab}.
\end{equation}
The integer $c_{ab}$ has a closed-form chord-geometric description.
List the distinct endpoints of $L_a \cup L_b$ in cyclic CCW order
starting from $L_a$'s first endpoint $p$, and write the resulting
arc edge-lengths between consecutive listed endpoints. Label each
endpoint by which chord(s) it belongs to: purely $L_a$, purely $L_b$,
or \emph{shared} (= in both, when the chords meet at a vertex).
Since $H = 2k+3$ is odd, the number of odd arcs is always odd (1 or 3).
The rule is:
\begin{itemize}
  \item If exactly one arc is odd and both its endpoints are strictly
    only-$L_a$ and only-$L_b$ respectively (no shared vertex involved),
    then $c_{ab} = -1$ when the arc goes from only-$L_a$ to only-$L_b$
    (CCW), and $c_{ab} = +1$ when it goes from only-$L_b$ to only-$L_a$.
  \item Otherwise $c_{ab} = 0$ — the unique odd arc is adjacent to a
    shared vertex, or is an interior arc of either chord, or there are
    three odd arcs.
\end{itemize}

\paragraph{Canonical Basis.}
The canonical basis consists of monomials of $L$'s labelled by
non-intersecting diagonals:
\begin{equation}
  L_{\{a_1, \cdots, a_n\}} \;\equiv\;
    \fq^{-\sum_{i<j} c_{a_i a_j}}\, L_{a_1} \cdots L_{a_n}\, .
\end{equation}

\paragraph{(X): quantum Ptolemy / Pl\"ucker.}
If $L_{\overline{ac}}$ and $L_{\overline{bd}}$ cross, they are the two
diagonals of a quadrilateral with vertices $a < b < c < d$ cyclically
on $P$. Let $(A_0, A_1, A_2, A_3)$ be the edge-lengths of the four arcs
$a \to b$, $b \to c$, $c \to d$, $d \to a$ (CCW). In the canonical
basis, quantum Ptolemy reads
\begin{equation}
  L_{\overline{ac}} \cdot L_{\overline{bd}}
  \;=\; \fq^{\alpha}\, L_{\{\overline{ab},\; \overline{cd}\}}
     \;+\; \fq^{\beta}\, L_{\{\overline{ad},\; \overline{bc}\}},
  \qquad (\alpha, \beta) \in \{(1, 0),\; (0, -1)\}.
\end{equation}
Since $H = 2k{+}3$ is odd, the two opposite arc-pairs $\{A_0, A_2\}$
(along $L_{\overline{ab}}, L_{\overline{cd}}$) and $\{A_1, A_3\}$ (along
$L_{\overline{ad}}, L_{\overline{bc}}$) differ in their numbers of
odd-length arcs:
$(\alpha, \beta) = (1, 0)$ if $\{A_0, A_2\}$ has fewer odd arcs,
$(0, -1)$ otherwise. When an opposite-edge pair degenerates to two
edges the corresponding canonical-basis element collapses to $1$.

\subsubsection{The trace}
The $\rho^2$-twisted trace $\Tr: A_\fq \to \bZ((\fq))$ is obviously $\rho^2$-invariant. As
$\rho^{2k+3}=1$, it is also $\rho$-invariant.  

There are $k+1$ elementary traces
\begin{equation}
  T_0 = \Tr(1), \qquad T_a = \Tr(L_{a, 0}), \quad a = 1, \ldots, k,
\end{equation}
and we found an algorithm to reduce all other trace to these by a judicious application of 
the twisted trace relation. As in the case of the pentagon algebra, the resulting linear combinations
have remarkable properties when combined with the $O(\fq)$ constraints. 

\paragraph{Reduction algorithm.}
Given a canonical-basis element $L = L_{\{a_1, \ldots, a_n\}}$ with $n \ge 2$, the trace 
$\Tr L$ is simplified recursively by a simple algorithm which cycles a generator repeatedly until 
it intersects other generators: we accumulate factors of $\fq$ through 
\begin{align}
L \xrightarrow{\text{factor}} L_{\{a_1, \ldots, a_{n-1}\}} L_{a_n} &\xrightarrow{\text{cycle}}  \rho^2(L_{a_n}) L_{\{a_1, \ldots, a_{n-1}\}} \xrightarrow{\text{commute across}} L_{\{a_1, \ldots, a_{n-1}\}}  \rho^2(L_{a_n}) \to \cr &\xrightarrow{\text{cycle}}  \rho^4(L_{a_n})L_{\{a_1, \ldots, a_{n-1}\}}\xrightarrow{\text{commute across}} \cdots \, ,
\end{align}
stopping when for some $j$ we find that $\rho^{2j}(L_{a_n})$ hits some $L_{a_i}$ it cannot pass through, and applying quantum Ptolemy to them. The quantum Ptolemy 
relation produces in each summand at least one generator which is {\it shorter} than both $L_{a_n}$ and $L_{a_i}$, so the procedure ultimately stops 
to a linear combination of elementary traces $T_a(\fq) = \Tr L_{a;i}$, which are $i$-independent, with coefficients which are Laurent polynomials in $\fq$. 

Each elementary trace is known as a linear combination of $M(2, 2k{+}3)$ minimal-model characters:
\begin{align}
  T_0 \;&=\; \chi_1(\fq^2), \\
  T_a \;&=\; (-1)^{m+1}\, \fq^{-m}\, \bigl(\chi_m - \chi_{m+1}\bigr)(\fq^2),
        \qquad m = m(a),
\end{align}
where $m = m(a) \in \{1, \ldots, k\}$ is given by
\begin{equation}
  m(a) \;=\;
    \begin{cases}
      a/2, & a \text{ even}, \\
      k - (a-1)/2, & a \text{ odd}.
    \end{cases}
\end{equation}

\subsubsection{The minors miracle}
The reduction algorithm of the previous subsection expresses every
canonical-basis trace as a $\bZ[\fq, \fq^{-1}]$-linear combination of
$T_0, T_1, \ldots, T_k$.  Fix a chord $L_{2;0}$, and choose $k-1$
other chords $D_1 = L_{1;0}$ and $D_{a-1}=L_{a>2,0}$ which do not intersect it and have distinct length. Form $k$
families of traces
\begin{equation}
  X^{(0)}_a := L_{2;0}^{\,a},
  \qquad
  X^{(r)}_a := L_{2;0}^{\,a}\cdot D_r, \quad r \neq 2,
  \qquad a = 1, 2, 3, \ldots,
\end{equation}
run the reduction algorithm on each, and collect the resulting
coefficient vectors
\begin{equation}
  \Tr(X^{(r)}_a) \;=\; \sum_{s = 0}^{k} c^{(r)}_s(a;\fq)\, T_s.
\end{equation}
Empirically, after dividing each $c^{(r)}_s(a;\fq)$ by the most negative overall power of $\fq$, the normalized coefficients stabilize as
$a \to \infty$ to honest $\fq$-series $M^{(r)}_s(\fq)$. By construction, 
\begin{equation}
  \sum_{s = 0}^{k} M^{(r)}_s(\fq)\, T_s = 0 \, ,
\end{equation}
and these equations turn out to be linearly independent. Clearly, this is enough to fix the ratios $T_s/T_0$,
and thus the full twisted trace up to scale. Remarkably, we find a stronger result: 
\begin{equation}
  (-1)^{s}\, \det \mathcal M^{(\hat s)} \;=\; T_s,
  \qquad s = 0, 1, \ldots, k,
\end{equation}
where $\mathcal M^{(\hat s)}$ is the matrix of $M^{(r)}_t$ with column $s$ deleted. We do not understand why this should be the case. 

\bibliographystyle{JHEP}

\bibliography{mono}

\end{document}